\documentclass[11pt]{article}

\usepackage{psfig}

\begin{document}

\mathsurround=1pt
\def\R{{\bf R}}
\def\RN{{{\bf R}^N}}
\def\Rn{{{\bf R}^n}}
\def\KN{{K^N}}
\def\b{{\vrule height6pt width6pt depth2pt}}
\def\MKN{{M\times K^N}}
\def\MRN{{M\times {\bf R}^N}}
\def\ksi{{\xi}}
\def\cv{{\cal V}}
\def\cvv{{{\cal V}_1  }}
\def\D{{\cal D}}
\def\Rp{{{\bf R}^p}}
\def\RnN{{{\bf R}^{n+N}}}
\def\Rs{{{\bf R}^s}}
\def\Rnsp{{\bf R}^{n-s-p}}
\def\Dp{{D^p}}
\def\Lt{{L^{\prime}}}

\newtheorem{thm}{Theorem}[section]
\newtheorem{defin}[thm]{Definition}
\newtheorem{cor}[thm]{Corollary}
\newtheorem{lem}[thm]{Lemma}
\newtheorem{prop}[thm]{Proposition}
\newtheorem{rem}[thm]{Remark}
\newtheorem{exam}[thm]{Example}
\newtheorem{sublem}[thm]{Sublemma}

\title {{\bf SURGERY ON LAGRANGIAN AND LEGENDRIAN SINGULARITIES}}

\author {Mikhail Entov\\
\\
Department of Mathematics\\ 
Stanford University\\
Stanford, CA 94305\\ 
USA\\
\\
\\
{\sl Current address:}\\
Department of Mathematics\\
Tel Aviv University\\
Ramat Aviv, Tel Aviv 69978\\
Israel\\
\\
e-mail: entov@math.tau.ac.il\\
}

\date{\today}
\maketitle

\begin{abstract}

Let
$\pi : E\to M$
be a smooth fiber bundle whose total space is a symplectic manifold 
and whose fibers are Lagrangian. Let
$L$
be an embedded Lagrangian submanifold of
$E$.
In the paper we address the following question: how can one simplify the 
singularities of the  projection
$\pi: L\to M$
by a Hamiltonian isotopy of
$L$
inside 
$E$?
We give an answer in the case when
${\rm dim}\ L = 2$ 
and  both
$L$
and
$M$
are orientable.
A weaker version of the result is proved in the higher-dimensional case.
Similar results hold in the contact category.

As a corollary one gets  an answer to one of the questions 
of V.Ar\-nold about the four cusps
on the caustic in the case of the Lagrangian collapse. As another corollary
we disprove  
Y.Chekanov's conjecture about singularities of the Lagrangian projection 
of certain Lagrangian tori in
${\bf R}^4$.

\end{abstract}

\section{Introduction and statement of main results}
\label{section-intro}

\subsection{Statement of the problem and an informal statement of results}
\label{ss-state-prob}

Let us start with the Lagrangian case.
Suppose that 
$\pi : E\to M$
is a
{\it
Lagrangian fiber bundle, 
}
i.e.
$E$
is a symplectic manifold and the fibers of
$\pi$
are Lagrangian submanifolds of
$E$.
The restriction of 
$\pi$
on a Lagrangian submanifold  
$L\subset E$
is called 
{\it
the Lagrangian projection of
$L$.
}
The singularities of the Lagrangian projection 
$\pi: L\to M$
are called 
{\it
the Lagrangian singularities of
$L$.
}

Suppose that
$L\subset E$
is an embedded Lagrangian submanifold. 
We would like to deform 
$L$
in order to simplify the singularities of its Lagrangian projection
$\pi: L\to M$
as much as we can.
More precisely, we are interested in the following question: 

{\it
Can we deform the Lagrangian submanifold
$L\subset E$
by a Hamiltonian isotopy
so that the singularities of the Lagrangian projection of the 
new submanifold on
$M$
are at worse folds?
And in general, what singularities can be obtained in this way?
}

This problem is motivated in particular by the V.Arnold's question asking
whether one can kill the Lagrangian cusp-singularities by a Hamiltonian
isotopy in the case of some special Lagrangian cylinder in
${\bf R}^4$
(called ``Lagrangian collapse'' -- see 
\cite{Ar2} 
and
Example~\ref{exam-0.1}
below for a detailed discussion).
Also Y.Chekanov conjectured in 
\cite{Ch}
that for some special Lagrangian 
tori in
${\bf R}^4$
that he had constructed 
the Lagrangian cusp-singularities cannot be removed by a Hamiltonian isotopy.

The answer to the questions above given in this paper can be formulated as 
follows. All the Lagrangian singularities higher
than folds (like cusps, swallowtails etc.) can be killed by a 
$C^0{\hbox{\rm -small}}$
Hamiltonian isotopy
{\it
maybe at the expense of creating additional folds,
}
provided that there is no obstruction for such a homotopy on the level
of 1-jets of the involved maps (in other words, if there is no
{\it
formal obstruction
}). 
We will prove this statement 
in the case when 
$L$
is a closed orientable two-dimensional surface and a weaker version of the
statement will be proved in the higher-dimensional case. 

The Legendrian case is completely parallel to the Lagrangian one (see
Section~\ref{ss-leg-case}).

Thus, Lagrangian (Legendrian) singularities higher than folds are
``soft'' from the point of view of symplectic (contact) topology -- one can 
get rid of them, if there is no formal obstruction for that. However 
the Lagrangian and Legendrian folds are ``rigid'' from the same point of 
view: they may not be killed, even if there is no formal obstruction for that.
For instance, for any
$n$
one can construct an embedded Legendrian sphere
$L\subset J^1 S^n$
which is smoothly isotopic and Legendrian regularly homotopic to the 
zero-section (and therefore there is no formal obstruction to kill all the folds by a Legendrian isotopy) but which is 
{\it
not
}
Legendrian isotopic to the zero-section
(see
\cite{En-lf}).
Legendrian submanifolds with analogous properties can also be constructed in
the space
$ST^\ast {\bf R}^n$
(see
\cite{En-lf}).
Similar examples have been independently constructed by 
P.Pushka${\rm r}'$ 
and E.Ferrand
(see 
\cite{Fer-Push}).

\bigskip
\bigskip
\noindent
{\bf An important assumption}

\bigskip
\noindent
{\sl
From the very beginning we assume that for all Lagrangian (Legendrian) 
submanifolds in this paper and all Lagrangian (Legendrian) singularities
or of singularities of tangency to a foliation, are contained in a compact 
part of the submanifold away from the boundary, if there is one. All 
deformations are supposed to be compactly supported and 
identical near the boundary. 
}
\bigskip
\bigskip

The singularities of the map
$\pi: L\to M$
can be interpreted as singularities of tangency of the submanifold
$L\subset E$
to the foliation 
$Vert$
of
$E$
by the fibers of 
$\pi:E\to M$.
According to the Lagrangian neighborhood theorem 
\cite{We1},\cite{W},
we can find a neighborhood
${\cal U}$
of the zero-section in
$T^\ast L$
and a map
$F:{\cal U}\to E$
which sends the zero-section onto 
$L$
and which is a symplectic embedding of
${\cal U}$
onto a neighborhood of
$L$ 
in
$E$.
Now, if we identify
$L$
with the zero-section in
$T^\ast L$,
then instead of singularities of tangency of
$L\subset E$
to
$Vert$
we can consider the singularities of tangency of the zero-section
$L\subset T^\ast L$
to the induced Lagrangian foliation
${\cal F} = F^\ast (Vert)$
of
${\cal U}$.

Thus instead of the original question we may consider the following
 harder problem:

\bigskip
\noindent
{\it
Let
${\cal F}$
be a Lagrangian foliation of a neighborhood 
${\cal U}$
of the zero-section 
$L$
in
$T^\ast L$.
How can the singularities of tangency of 
$L$
to
${\cal F}$
change when we deform
$L$
by a Hamiltonian isotopy inside 
${\cal U}$?
}
\bigskip

\subsection{The structure of Lagrangian singularities: basic definitions}
\label{ss-str-lagr-sing}

First we introduce a notational convention:
by a 
{\it
Lagrangian distribution
} 
along an embedded Lagrangian submanifold
$\Lambda\subset T^\ast L$
we always mean a field of Lagrangian planes in
$T (T^\ast L)$
defined on
$\Lambda$. 
If
${\cal N}$
is a Lagrangian distribution along
$\Lambda$, 
we write
${\cal N}_x$
for a fiber of
${\cal N}$
at
$x\in \Lambda$.

\begin{defin}
\label{def-0.1}
{\rm
We say that a Lagrangian distribution
$\cal N$ 
along
$L$
is
\break
{\it
$\Sigma^2{\hbox{\it -nonsingular,}}$
}
if at any point 
$x\in L$
we have
${\rm dim}\ (T_x L\cap {\cal N}_x)\leq 1$.
Otherwise
$\cal N$
is called
$\Sigma^2{\hbox{\it -singular}}$.
}
\end{defin}

The term
$\Sigma^2$
comes from the Thom-Boardman hierarchy of singularities
(see \cite{Th}, \cite{B}).
According to the singularity theory (see e.g.
\cite{AGV})
the singularities of tangency of a Lagrangian submanifold
$L\subset T^\ast L$
to a generic
$\Sigma^2{\hbox{\rm -nonsingular}}$
Lagrangian distribution 
${\cal N}$
form a ``complete flag'' 
$L^n = V_0\supset V_1\supset V_2\supset\ldots\supset V_{n+1}=\emptyset$
of (not necessarily connected) smooth submanifolds of
$L$.
Namely, along 
$L\setminus V_1$
the distribution 
${\cal N}$
is transversal to
$L$,
along
$V_1$
the distribution planes cut out a line field in the tangent planes of 
$L$,
and then
$V_2\subset V_1$
is exactly the set of points where this line field is not transversal to
$V_1$,
$V_3\subset V_2$
is exactly the set of points where the line field is not transversal to
$V_2$
etc.
According to the Thom-Boardman terminology, along 
$L\setminus V_1$
the singularities of tangency are of type 
$\Sigma^1$,
along 
$V_2\setminus V_3$
the singularities are of type 
$\Sigma^{11}$
etc.

For 
$i\geq 2$
there exists a natural coorientation of
$V_{i-1}$
in
$V_{i-2}$
along 
$V_i\setminus V_{i+1}$,
(see 
Definition~\ref{def-1.4}),
and the submanifold
$V_1$
has a natural coorientation in
$L$
coming from the
{\it 
Maslov coorientation
}
of Lagrangian folds
(see 
\cite{Ar4}).
We express all these coorientations by means of some unit normal vector fields
$v_i$
(with respect to some Riemannian metric fixed once and forever).
Thus for
$i\geq 2$
the unit vector field
$v_i$
is defined along
$V_i\setminus V_{i+1}$
and is normal to
$V_{i-1}$
in
$V_{i-2}$,
and the unit vector field
$v_1$
is defined along
$V_1$
and is normal to
$V_1$
in
$L$.

\bigskip
\begin{defin}
\label{def-0.2}
{\rm
A
$\Sigma^2{\hbox{\rm -nonsingular}}$
Lagrangian dis\-tri\-bu\-tion
$\cal N$ 
on
$L$
is called a
{\it
fold-type Lagrangian distribution,
}
if the submanifold
$V_1$
in the flag above is smooth and
${\rm dim}\ ({\cal N}_x \cap T_x V_1) = 0$
for any
$x\in V_1$.
In other words, 
$V_i=\emptyset$
for any
$i>1$.
}
\end{defin}
\bigskip

Suppose that we have a 
$\Sigma^2{\hbox{\rm -nonsingular}}$
Lagrangian distributions
${\cal N}$
and
${\cal N}^\prime$
along, respectively, some Lag\-ran\-gian submanifolds
$L$
and
$L^\prime$
in
$T^\ast M$.
Let
$\{V_i,v_i\}$
(resp.
$\{V_i^\prime, v_i^\prime\}$)
denote the data describing as above the singularities of tangency of
${\cal N}$ 
to
$L$
(resp.
${\cal N}^\prime$ 
to
$L^\prime$).

We say that the singularities of tangency 
of 
${\cal N}$
to
$L$
are 
{\it
equivalent 
}
to the singularities of tangency of 
${\cal N}^\prime$
to
$L^\prime$,
if there exists an isotopy between
$L$
and
$L^\prime$
in
$T^\ast M$
which maps the data
$\{V_i,v_i\}$
into the corresponding data
$\{V_i^\prime, v_i^\prime\}$
(also see Definition~\ref{def-1.3}).

\subsection{Double folds}
\label{ss-df}

Before stating the main results of the paper 
we need to define a particular combination of singularities.

\begin{defin}
\label{double-folds}
{\rm
Let
${\cal N}$
be a Lagrangian distribution along 
$L$.
Let
$S_1^{n-1}$, 
$S_2^{n-1}\subset L$
be embedded spheres such that
$S_2^{n-1}$
bounds a ball 
$B^n$
in
$L$
and the sphere
$S_1^{n-1}$
lies inside that ball.

We say that
$L$
has a 
{\it
double fold 
}
with respect to 
${\cal N}$
along the spheres
$S_1^{n-1}, S_2^{n-1}\subset L$,
if the restriction of
${\cal N}$
on a neighborhood of the ball
$B^n$
is a fold-type distribution that has singularities of tangency to
$L$
only along the spheres
$S_1^{n-1}$
and
$S_2^{n-1}$
and the (Maslov) coorientations of both spheres in
$L$
point in the opposite directions.

In general, a pair of spheres
$S_1^{n-1}, S_2^{n-1} \subset L$
as above equipped with opposite coorientations will be called a
{\it
double fold pair of (cooriented) spheres.
}
}
\end{defin}

\begin{figure}
\centerline{\psfig{figure=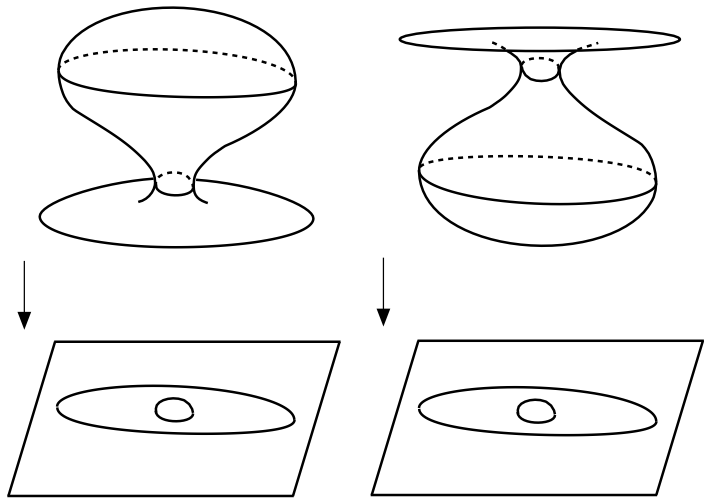,height=80mm}}
\caption{Double folds.}
\label{fig11}
\end{figure}

\begin{rem}
\label{d-fold-n-form}
{\rm
In a small neighborhood 
$U$
of the ball
$B^n$
from Definition~\ref{double-folds}
we can reduce 
${\cal N}$
and 
$L$
by a symplectomorphism to a specific form so that

\smallskip
\noindent
(i) the neighborhood
$U$
is presented as an open set in the standard symplectic space
${\bf R}^{2n}$
with the standard Darboux coordinates
$p_1,\ldots ,p_n , q_1 ,\ldots ,q_n$;

\noindent
(ii) the Lagrangian distribution
${\cal N}$
is identified with a field of affine Lagrangian planes parallel to the 
vertical 
$(p_1,\ldots ,p_n){\hbox{\rm -plane}}$;

\noindent
(iii) the Lagrangian submanifold
$L$
has a  ``mushroom''-type shape with respect to the vertical plane field 
${\cal N}$
as on 
Fig.~\ref{fig11}.

\noindent
See 
Section~\ref{double-folds-pfs} 
for details.
}
\end{rem}
\bigskip

\subsection{The main results}
\label{ss-main-res}

Let us identify 
$L$
with
the zero-section in
$T^\ast L$
and let
${\cal F}$
be a Lagrangian foliation defined on a neighborhood of
$L$
in
$T^\ast L$.
Also let
${\cal N}$
be a generic
$\Sigma^2{\hbox{\rm -nonsingular}}$
Lagrangian distribution along
$L$.

We first 
state the result for
the case when
${\rm dim}\ L = 2$
and all manifolds and distributions are orientable (Theorem A). 
Then we present a weaker version of the result for a general
$L$
of any dimension (Theorem B).

Lagrangian singularities of the Thom-Boardman type
$\Sigma^2$
generically have codimension 3
(see 
\cite{Ar1},\cite{AGV}). 
Therefore a generic Lagrangian distribution along a
two-dimensional surface is
$\Sigma^2{\hbox{\rm -nonsingular}}$.

\bigskip
\noindent
{\bf Theorem A}
{\it
Suppose that
$L$
is an oriented surface and that the Lagrangian foliation
${\cal F}$
and the Lagrangian distribution 
${\cal N}$
are oriented. 
Suppose also that the oriented Lagrangian distributions
${\left. {\cal F}\right|}_L$ 
and
${\cal N}$
along 
$L$
can be connected in the class of  such distributions.

Then there exists a Hamiltonian isotopy
$\{h_t (L)\}$,
$0\leq t \leq 1$,
such that the singularities of tangency of the Lagrangian submanifold 
$h_1(L)\subset T^\ast L$
to
${\cal F}$
are equivalent to the union of the singularities of tangency of
$L$
to
${\cal N}$
and a number of additional double fold pairs of circles in
$L$.

The Hamiltonian isotopy above can be chosen in such a way that it is 
compactly supported and arbitrarily
$C^0{\hbox{\it -small}}$
and so that 
$\displaystyle \max_{t,\gamma}
\mid \int_{h_t(\gamma)} pdq \mid$ 
is also arbitrarily
$C^0{\hbox{\it -small}}$
(by
$\gamma$
we denote a path in
$L$).
}
\bigskip

\begin{rem}
\label{rem-01}
{\rm
The hypothesis of Theorem A is a necessary condition for the existence of a 
$C^0{\hbox{\rm -small}}$
Hamiltonian isotopy which deforms
$L$
into a Lagrangian submanifold whose singularities of tangency to
${\cal F}$
are equivalent to the singularities of tangency of
$L$
to
${\cal N}$.
Indeed, if there exists such a Hamiltonian isotopy 
$\{ h_t (L)\}$,
$0\leq t\leq 1$,
then there exists a homotopy 
$\{ {(h_t^{-1})}_\ast ({\cal F})\}$
in the class of Lagrangian distributions along 
$L$
connecting 
${\left. {\cal F}\right|}_L$
with a Lagrangian distribution
${(h_1^{-1})}_\ast ({\cal F})$,
so that the singularities of tangency of
${(h_1^{-1})}_\ast ({\cal F})$
to
$L$
are equivalent to the singularities of tangency of
${\cal N}$
to
$L$.

}
\end{rem}
\bigskip

In the higher-dimensional case one has the following weaker result 
(as compared to Theorem A).

\bigskip
\noindent
{\bf Theorem B}
{\it
Suppose that
${\rm dim}\ L\geq 2$.
Suppose also that the Lagrangian distribution
${\left. {\cal F}\right|}_L$
along
$L$
is generic and
$\Sigma^2{\hbox{\it -nonsingular}}$ 
and that it can be connected with
${\cal N}$
in the class of 
$\Sigma^2{\hbox{\it -nonsingular}}$ 
Lagrangian distributions along
$L$.

Then there exists a Hamiltonian isotopy
$\{h_t (L)\}$,
$0\leq t \leq 1$,
such that the singularities of tangency of the Lagrangian submanifold 
$h_1(L)\subset T^\ast L$
to
${\cal F}$
are equivalent to the union of the singularities of tangency of
$L$
to
${\cal N}$
and a number of additional double fold pairs of spheres in
$L$.

In particular, if
${\cal N}$
is a fold-type distribution, then
all the singularities of tangency of
$h_1(L)$
to
${\cal F}$
are folds.

The Hamiltonian isotopy above can be chosen in such a way that it is 
compactly supported and arbitrarily
$C^0{\hbox{\it -small}}$
and so that 
$\displaystyle \max_{t,\gamma}
\mid \int_{h_t(\gamma)} pdq \mid$ 
is also arbitrarily
$C^0{\hbox{\it -small}}$.
}
\bigskip

Unlike in Theorem A, the hypothesis of Theorem B is not a necessary condition 
for the existence of a Hamiltonian isotopy that deforms 
$L$
into a Lagrangian submanifold whose singularities of tangency to
${\cal F}$
are equivalent to the singularities of tangency of
$L$
to
${\cal N}$.
The reason is that we require 
the homotopy between 
${\left. {\cal F}\right|}_L$
and
${\cal N}$
to lie in the class of 
$\Sigma^2{\hbox{\rm -nonsingular}}$ 
and not in the class of
{\sl 
all
}
Lagrangian distributions along 
$L$
(see Remark~\ref{rem-01}).

\subsection{Corollaries of Theorem A}
\label{ss-cor-2dim}

First we explicitly describe the possible sets of Lagrangian 
singularities on a given closed Lagrangian surface.

\begin{figure}
\centerline{\psfig{figure=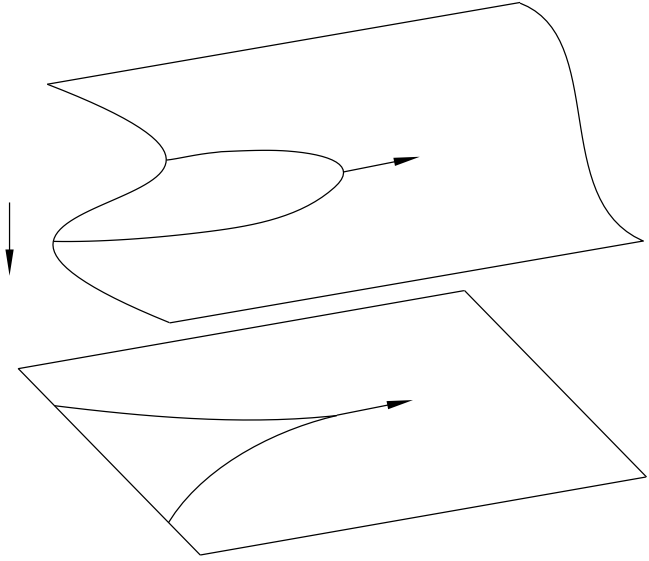,height=80mm}}
\caption{The characteristic direction of a cusp on the caustic and the corresponding vector at the cusp-singularity point of the Lagrangian submanifold.}
\label{fig2}
\end{figure}

If
$L$
is generic,
its Lagrangian singularities form a closed 1-dimensional
submanifold
$V_1\subset L$
(the closure of the set of Lagrangian folds)
with a 0-dimensional submanifold
$V_2\subset V_1$
of Lagrangian cusps.
The submanifold 
$V_1$
is equipped with a coorienting unit vector field 
$v_1$,
corresponding to the Maslov coorientation of the Lagrangian folds in
$L$.
The points of the zero-dimensional submanifold
$V_2$ 
are equipped with unit vectors 
(denoted altogether by
$v_2$)
normal to
$V_1$, 
that arise from the characteristic directions of cusps on the caustic  
(i.e. on the set of critical values of the Lagrangian projection) -- see
Fig.\ref{fig2}.
Thus the structure of the Lagrangian singularities is described by the 
data
$\{ V_1, V_2, v_1,v_2\}$.

Since
$V_1$
is cooriented in 
$L$
it subdivides
$L$
into two parts
$L_1$
and
$L_2$
with the common boundary
$V_1$.
Suppose the
$V_2$ 
consist of
$n_1 + n_2$
points so that at
$n_1$
of the points the vectors 
$v_2$ 
are directed into 
$L_2$
and at the rest
$n_2$
points into 
$L_1$.

Now suppose that
$L$
is an embedded closed oriented Lagrangian submanifold in
$T^\ast M$, 
where
$M$
is an oriented surface.
Let
$c (L)\in H^1(L)$
be the Maslov class and let
$d$
be the degree of the Lagrangian projection
$L\to M$
(in particular, if 
$M$
is not closed,
$d=0$).
It is well known that

\smallskip
\noindent
(A) $V_1$
realizes a homology class Poincar{\'e}-dual to
$c(L)$
($V_1$
is cooriented and therefore also oriented because
$L$
is oriented);

\noindent
(B) 
$\chi (L_1) - \chi (L_2) + n_1 - n_2 = d\chi (M)$
(this formula holds for any smooth map
$L\to M$,
see e.g. 
\cite{Ha},\cite{El2}).

\begin{cor}
\label{cor-0.11}
Suppose that a formal data 
$\{ V_1, V_2, v_1,v_2\}$
on 
$L$
satisfies the conditions (A) and (B) above.
Then there exists a (generic) Lagrangian submanifold
$L^\prime\subset T^\ast M$
Hamiltonian isotopic to
$L$
such that the singularities of the Lagrangian projection
$L^\prime\to M$
are equivalent to the union of
$\{ V_1, V_2, v_1,v_2\}$
and of a number of double fold pairs of cooriented circles in
$L$.

The Hamiltonian isotopy above can be chosen in such a way that it is 
compactly supported and arbitrarily
$C^0{\hbox{\it -small}}$
and so that   
  $\displaystyle \max_{t,\gamma}
\mid \int_{h_t(\gamma)} pdq  - \int_{\gamma} pdq  \mid$ 
is also arbitrarily
$C^0{\hbox{\it -small}}$.
\end{cor}

\bigskip
\begin{cor}
\label{cor-0.1}
Let
$L$
be an embedded orientable Lagrangian submanifold in
$T^\ast M$, 
where
$M$
is an orientable surface.
Then 
$L$
is Hamiltonian isotopic to a Lagrangian submanifold
$L_1\subset T^\ast M$
such that all the singularities of the Lagrangian projection
$L_1\to M$
are folds.

The Hamiltonian isotopy can be chosen in such a way that it is 
compactly supported and arbitrarily
$C^0{\hbox{\it -small}}$
and so that   
  $\displaystyle \max_{t,\gamma}
\mid \int_{h_t(\gamma)} pdq  - \int_{\gamma} pdq  \mid$ 
is also arbitrarily
$C^0{\hbox{\it -small}}$.

\end{cor}
\bigskip

Corollaries~\ref{cor-0.11},\ref{cor-0.1} 
follow from Theorem A and from the general results 
(see Propositions~\ref{prop-22.1},\ref{prop-22.2},\ref{prop-22.3}) 
which classify the homotopy classes of oriented Lagrangian 
distributions along an oriented surface and describe what 
singularities can be realized as singularities of tangency of
$L$
to an oriented 
$\Sigma^2{\hbox{\rm -nonsingular}}$
Lagrangian distribution
(see Section~\ref{two-dim-case}).

Corollary~\ref{cor-0.1},
in particular, disproves a conjecture by 
Y.Che\-ka\-nov 
\cite{Ch}
that for certain Lagrangian tori in
${\bf R}^4$
the Lagrangian projection necessarily has at least 4 cusp singularities.

Now let us consider the case when 
$i:L\to T^\ast M$\
is a Lagrangian immersion
and both
$L$
and
$M$
are orientable surfaces.
It follows from the Lagrangian neighborhood theorem
\cite{We1},\cite{W}
that there exists a symplectic immersion 
$F$
of a neighborhood
${\cal U}$
of
$L$
in
$T^\ast L$
onto a neighborhood of
$i(L)$
in
$T^\ast M$.
Then one has the following 
version of Corollary~\ref{cor-0.1}.

\begin{cor}
\label{cor-0.2}
Let
$i:L\to T^\ast M$ 
and
$F: {\cal U}\to T^\ast M$
be as above.
Then there exists a  Hamiltonian isotopy of
$L$
into
$L_1$
inside
${\cal U}\subset T^\ast L$
so that all singularities of the Lagrangian projection of
an immersed Lagrangian submanifold
$F(L_1)\subset T^\ast M$
on
$M$
are folds.

The Hamiltonian isotopy can be chosen in such a way that it is 
compactly supported and arbitrarily
$C^0{\hbox{\it -small}}$
and so that 
$\displaystyle \max_{t,\gamma}
\mid \int_{F\circ h_t(\gamma)} pdq  - \int_{F(\gamma)} pdq  \mid$ 
is also arbitrarily
$C^0{\hbox{\it -small}}$.
\end{cor}
\bigskip

\subsection{Importance of some of the assumptions in the results above}

\begin{rem}
\label{rem-0.1}
{\rm
The claims of Corollaries~\ref{cor-0.1},\ref{cor-0.2}
are not necessarily true, if we do not require the orientability of both
$L$
and
$M$.
For instance, if 
$M$ 
is an orientable surface without boundary,
$L$
is closed and 
$\chi (L)$
is odd (which implies that
$L$ 
is  non-orientable), then
a purely topological obstruction dictates that the number 
of cusp-singularities for
{\it 
any smooth map 
}
from
$L$
to
$M$
is necessarily odd (see
\cite{Wh},\cite{Th},\cite{Ha},\cite{Le},\cite{El2}). 
}
\end{rem}
 
\begin{rem}
\label{rem-02}
{\rm
The claims of Theorems A,B and of their 
Corollaries~\ref{cor-0.11}, \ref{cor-0.1}, \ref{cor-0.2} 
might not be true, if one does not 
allow additional double folds -- see 
the discussion about the ``rigidity'' of folds in
Section~\ref{ss-state-prob}.

}
\end{rem}
\bigskip

\subsection{The case of a Lagrangian manifold defined by a 
generating function. Lagrangian collapse}
\label{ss-gfcn}

We recall the definition of 
{\it 
a generating function
}
for a Lagrangian (Legendrian) submanifold.

\begin{defin}
\label{gen-fcns-defin}
{\rm
Let
$K$
be a manifold and let
$f:M\times K\to {\bf R}$
be a function.
Consider the graph 
${\rm Graph} (df)$
of 
$df$
in
$T^\ast (M\times K)$.
This is an exact Lagrangian submanifold 
of the symplectic space
$T^\ast (M\times K)$.
Consider a coisotropic surface  
$\Sigma\subset T^\ast (M\times K)$
defined as the set of those 1-forms on
$M\times K$
which vanish along the fibers of the projection
$M\times K\to K$.
The characteristic foliation on
$\Sigma$
provides a fibration of 
$\Sigma$
over 
$T^\ast M$.

We apply to
${\rm Graph} (df)$
and
$\Sigma$
the standard construction called the 
{\sl
Lagrangian reduction}.
Namely, suppose that  
${\rm Graph} (df)$
intersects 
$\Sigma$
transversally.
Then the intersection 
${\rm Graph} (df)\cap\Sigma$
is a smooth manifold which under the fibration
$\Sigma\to T^\ast M$
projects into an immersed exact Lagrangian submanifold
$L\subset T^\ast M$.
In such a situation
$f$
is called a 
{\it
generating function}
for
$L$.

If instead of 
${\rm Graph} (df)$
we take a graph of a closed 1-form
$\lambda$
on
$M\times K$,
such that the graph intersects 
$\Sigma$
transversally, 
then the same construction gives us an immersed 
(and not necessarily exact) Lagrangian submanifold of
$T^\ast M$, 
which, we say, is 
{\it
defined by the 
generating form
$\lambda$}.

Using a similar construction (called the 
{\sl
Legendrian reduction
})
one can start with a function on
$M\times K$,
consider a Legendrian submanifold in the 1-jet space
$J^1 (M\times K)$
which is the graph of the 1-jet of
$f$
and define a 
{\it
a Legendrian submanifold of 
$J^1 M$
generated by the function 
$f:M\times K\to {\bf R}$
}.

The Lagrangian submanifold in
$T^\ast M$
generated by a function
$f:M\times K\to {\bf R}$
is defined in local coordinates by the formula

\[ 
\left\{ 
\begin{array}{cc}
p = {\partial f}(q,\xi)/{\partial q}\\
0 = {\partial f}(q,\xi)/{\partial \xi}
\end{array}
\right.
\]

A Legendrian submanifold in
$J^1 M$
generated by 
$f:M\times K\to {\bf R}$
is defined in local coordinates by the formula

\[ 
\left\{ 
\begin{array}{lll}
z = f(q,\xi)\\
p = {\partial f}(q,\xi)/{\partial q}\\
0 = {\partial f}(q,\xi)/{\partial \xi}
\end{array}
\right.
\]

Here 
$q,p$
(resp.
$q,p,z$)
are the standard local Darboux coordinates on
$T^\ast M$
(resp.
$J^1 M$)
and
$\xi$
stands for some local coordinates on
$K$.

}
\end{defin}
\bigskip

The following important example is 
due to V.Arnold, who calls it  ``Lagrangian collapse'' 
\cite{Ar2}.

\begin{exam}
\label{exam-0.1}
{\rm
Let
$L$
be an embedded Lagrangian submanifold of 
$T^\ast {\R}^2$
generated by the function
$f(q_1 , q_2 ,\xi)=q_1 cos \xi +q_2 sin \xi$
defined on
${\R}^2\times S^1$.
The question posed by V.Arnold (\cite{Ar2}, also see \cite{Ar5}) is whether 
the caustic of Lagrangian projection of a (generic)
Lagrangian submanifold Hamiltonian isotopic to
$L$ 
in
$T^\ast {\R}^2$
always has at least 4 cusps.
In  
\cite{Ar2}
V.Arnold proved that there must be at least 4 such cusps, if the
Hamiltonian isotopy is sufficiently
$C^1{\hbox{\rm -small}}$.
Corollary~\ref{cor-0.1}
gives a negative answer to the question in the general case.
In Section~\ref{sect-gen-f-pfs}
we present an explicit sequence of surgeries on
the singularities of the Lagrangian projection of 
$L$
that eliminate all the cusps on the caustic (and each of these 
surgeries can be realized by  a 
$C^0{\hbox{\rm -small}}$
Hamiltonian isotopy of
$L$).
S.Ivankov previously constructed a deformation of
$L$ 
in the class of
{\it immersed}
Lagrangian submanifolds which eliminates all the cusps on the caustic 
(unpublished).

}
\end{exam}
\bigskip

In conjunction with Example~\ref{exam-0.1} 
one might also ask  
the following question:  is it possible to choose the 
$C^0{\hbox{\rm -small}}$
Hamiltonian isotopy destroying the cusps 
in 
Example~\ref{exam-0.1} 
so that it is covered by a deformation of the generating function
$f$
in the class of generating functions defined on
${\R}^2\times S^1$?
The following proposition provides an answer to this question.

\begin{prop}
\label{thmC}
Suppose that an embedded
Lagrangian submanifold
$L\subset T^\ast M$
is generated by a function
$f$
defined (on an open subset of) 
$M\times K^1$,
where
$K$
is a (connected)
$1{\hbox{\it -dimensional}}$
manifold.

Then 
$L$
can be deformed by a Hamiltonian isotopy into a new Lagrangian submanifold 
$L_1$
such that the canonical Lagrangian projection
of
$L_1$
on
$M$
has only fold-type singularities. Moreover,
this Hamiltonian isotopy can be covered by a 
deformation 
$\{ f_t\}$
of the generating function
$f$ 
in the class of generating functions
defined (on an open subset of)
$M\times K^1$.

The deformation 
$\{ f_t\}$
can be made compactly supported and arbitrarily
$C^1{\hbox{\it -small}}$.
\end{prop}

\begin{rem}
\label{rem-0.100}
{\rm
Statements similar to Proposition~\ref{thmC} 
can be proved (by exactly the same method as 
Proposition~\ref{thmC})
for a Lagrangian submanifold generated by a 1-form defined on
$M\times K^1$
and also in the case of an Lagrangian submanifold.

}
\end{rem}

\begin{cor}
The 
$C^0{\hbox{\it -small}}$
Hamiltonian isotopy destroying the cusps in 
Example~\ref{exam-0.1}
can be chosen so that it is covered by a 
deformation of the generating function 
$f(q_1 , q_2 ,\xi)=q_1 cos \xi +q_2 sin \xi$
as a function defined on
${\R}^2\times S^1$
(see 
Example~\ref{exam-2.1}
in 
Section~\ref{sect-gen-f-pfs}).
The deformation of the function
$f$
is compactly supported and can be made arbitrarily
$C^1{\hbox{\it -small}}$
\end{cor}

\begin{rem}
\label{rem-0.200}
{\rm
One can prove a strengthened version of 
Proposition~\ref{thmC}  
similar to Theorem B 
that says  what kind of singularities
can be realized (up to possible additional double folds)  by a deformation of
a generating function
$f$ 
in the class of generating functions 
defined on (an open subset of)
$M\times K^1$ --
see Proposition~\ref{thmC-full-str}.
}
\end{rem}

\subsection{The Legendrian case}
\label{ss-leg-case}

Now let us discuss the Legendrian case. 
Suppose that 
$\pi : E\to M$
is a
{\it
Legendrian fiber bundle, 
}
i.e.
$E$
is a contact manifold and the fibers of
$\pi$
are Legendrian submanifolds of
$E$.
The restriction of 
$\pi$
on a Legendrian submanifold  
$L\subset E$
is called 
{\it
the Legendrian projection of
$L$.
}
The singularities of the Legendrian projection 
$\pi: L\to M$
are called 
{\it
the Legendrian singularities of
$L$.
}

Suppose that
$L\subset E$
is an embedded Legendrian submanifold. 
In the same fashion as before
we can use the Legendrian neighborhood theorem 
\cite{Ly}
(similar to the Lagrangian neighborhood theorem of Weinstein)
and reduce the question about the singularities of a Legendrian projection
$L\to M$ 
to a question about the singularities of tangency of the zero-section in
$J^1 L$
to a Legendrian foliation defined on a neighborhood of the zero-section.

Observe that under the canonical projection
$J^1 L\to T^\ast L$
Legendrian planes are mapped into Lagrangian planes
and this projection establishes a one-to-one
correspondence between 
Legendrian submanifolds 
of
$J^1 L$
and exact immersed Lagrangian submanifolds of
$T^\ast L$.
The singularities of tangency of a Legendrian submanifold of
$J^1 L$
to a Legendrian distribution are the same as for the corresponding
(immersed) Lagrangian submanifold and the corresponding Lagrangian 
distribution downstairs in
$T^\ast L$.
(A Legendrian distribution along a submanifold
$\Lambda\subset J^1 L$
is a field of Legendrian planes in
$T (J^1 L)$
defined on
$\Lambda$.)
A 
$C^0{\hbox{\rm -small}}$
Hamiltonian isotopy 
$\{ h_t (L)\}$
of the zero-section of
$T^\ast L$,
such that 
$\displaystyle \max_{t,\gamma}
\mid \int_{h_t(\gamma)} pdq\mid$ 
is small, corresponds to a
$C^0{\hbox{\rm -small}}$
contact isotopy of the zero-section of
$J^1 L$. 
In general, we can carry over everything we said before for 
the Lagrangian case to the Legendrian
case by projecting a Legendrian submanifold
$\Lambda$
into
$T^\ast L$,
considering the corresponding Lagrangian (symplectic) 
object there and then lifting it upstairs to
$J^1 L$
as a Legendrian (contact) object.

{\sl
Thus all theorems, corollaries and remarks above for the Lagrangian 
case have direct analogs in the Legendrian case.
}

For example, the following Legendrian analog of Theorem A can 
be immediately deduced from Theorem A. Let us identify
$L$
with the zero-section of
$J^1 L$.
Let
${\cal F}$
be a Legendrian foliation defined on a neighborhood of
$L$
in
$J^1 L$,
and let
${\cal N}$
be a generic 
$\Sigma^2{\hbox{\rm -nonsingular}}$
Legendrian distribution along
$L$
(the definition of a
$\Sigma^2{\hbox{\rm -nonsingular}}$
Legendrian distribution can be given in the same fashion as 
Definition~\ref{def-0.1}
in the Lagrangian case).

\bigskip
\noindent
{\bf Theorem A'}
{\it
Let
$L$
be an oriented surface and suppose that the Legendrian foliation
${\cal F}$
and the Legendrian distribution
${\cal N}$
are oriented. 
Suppose also that the oriented Legendrian distributions
${\left. {\cal F}\right|}_L$ 
and
${\cal N}$
along 
$L$
can be connected in the class of such distributions.

Then
$L$
is contact isotopic to
a Legendrian submanifold 
$L_1\subset J^1 L$
such that the singularities of tangency of
$L_1$
to
${\cal F}$
are equivalent to the union of the singularities of tangency of
$L$
to
${\cal N}$
and a number of additional double fold pairs of spheres in
$L$.

The contact isotopy can be made compactly supported and arbitrarily
$C^0{\hbox{\it -small}}$.

}
\bigskip

The Legendrian counterpart for Corollary~\ref{cor-0.1}
has been independently announced by P.Pu\-shka${\rm r}'$.

\subsection{The sketch of basic ideas used in the proofs of the main results.}
\label{ss-basic-ideas}

The methods we use to prove theorems originate from
\cite{El}. 
First, we prove Theorem B using the following idea. As we said in
Section~\ref{ss-str-lagr-sing},
a description of singularities of tangency of 
$L$
to a generic
$\Sigma^2{\hbox{\rm -nonsingular}}$
Lagrangian distribution along
$L$
can be given by a ``complete flag'' 
$L^n\supset V_1\supset V_2\supset\ldots\supset V_{n+1}=\emptyset$
of (not necessarily connected) smooth submanifolds
of
$L$,
equipped with some vector fields
$v_i$.
In general, such an object (a ``complete flag'' of submanifolds equipped 
with appropriate vector fields) will be called a
{\it
chain
}
(see 
Section~\ref{ch-dirsurg}).

Consider a homotopy 
$\{ {\cal N}^t\}$
connecting
${\cal F}$
with 
${\cal N}$
in the class of 
\break
$\Sigma^2{\hbox{\rm -nonsingular}}$  
Lagrangian distributions along
$L$
(see Theorem B). 
We follow the homotopy 
$\{ {\cal N}^t\}$
and see what kind of bifurcations happen to the chains associated with each
${\cal N}^t$.
Since we are working only with 
\break
$\Sigma^2{\hbox{\rm -nonsingular}}$  
Lagrangian distributions, we can use the classification of the
{\it typical}
bifurcations of the singularities of tangency of
$L$
to
${\cal N}^t$, 
i.e. the classification of 
the typical surgeries on the submanifolds from a chain
associated with each Lagrangian distribution
${\cal N}^t$
(see 
\cite{Ar3},\cite{AGV}).
These surgeries can be, roughly speaking, 
of ``birth''
and ``death'' types (following the terminology from \cite{El} we call them 
{\it direct
}
and 
{\it 
inverse 
}
surgeries). 
By possibly 
creating additional spheres of folds we can always
decompose a ``death''-type bifurcation into 
a product of ``birth''-type ones and, thus from the homotopy
$\{ {\cal N}^t\}$
we obtain a ``code'' -- a sequence of
``birth''-type bifurcations 
(direct surgeries)  -- see 
Section~\ref{ch-dirsurg}.

Then we consider the singularities of tangency of
${\cal F}$
to
$L$
and, given a ``code'', or a sequence of bifurcations 
(direct surgeries), we deform 
$L$
in
$T^\ast L$
in several steps so that the singularities of tangency of
${\cal F}$
to
$L$
change at each step according to the corresponding
``prescribed'' bifurcation (direct surgery) from the ``code''.

Proposition~\ref{thmC} 
is proved by similar methods. To prove Theorem A we 
reduce it by certain tricks to Theorem B.

\subsection{A brief review of the paper by sections}.
\label{ss-paper-org}

The paper is organized as follows.
In 
Section~\ref{ch-dirsurg} 
we introduce the notions of a chain and of a direct surgery which 
are our main tools in dealing with singularities of maps. In 
Section~\ref{sect-pf-thmB} 
we prove 
Theorem B 
modulo the technical details which are discussed in 
Section~\ref{pf-prop-3.1-chap}. 
In 
Section~\ref{two-dim-case}
we discuss the two-dimensional case and prove Theorem A and the Corollaries~\ref{cor-0.11},\ref{cor-0.1}. 
In 
Section~\ref{gen-f} 
we introduce  the necessary setup for the case when a Lagrangian submanifold 
is defined by a generating function. In 
Section~\ref{sect-gen-f-pfs} 
we prove 
Proposition~\ref{thmC} 
and discuss in details the example of Lagrangian collapse. In 
Section~\ref{pf-prop-3.1-chap}
we prove our main technical result which shows that a direct surgery can 
be realized by a deformation of a generating function. In 
Section~\ref{double-folds-pfs} 
we explain how to create double folds and what are the ``normal forms'' 
for them.

\subsection{Acknowledgments} 
\label{ss-ackn}

I am deeply grateful to Y.Eliashberg for the
numerous helpful discussions and encouragement. 
He suggested the problem to me and generously and patiently shared with me his
ideas on the subject. Without him this work would have never been completed.
I am grateful to M.Kazarian for pointing out  the mistakes in
an earlier version of the paper and for his help in correcting them. 
I thank Y.Chekanov, N.Mishachev, D.Novikov and V.Zakalyukin 
for the useful conversations and encouragement. I am indebted to L.Polterovich
and to the referee for their remarks on the manuscript.

\section{Chains and direct surgeries}
\label{ch-dirsurg}

In this section we recall some basic notions concerning the structure
of singularities (see 
\cite{El}).
Let
$V^n$
be an (immersed) submanifold  of
$E^{n+N}$
$(n\geq 2)$. 
Let us assume for simplicity that the manifold
$E$
is equipped with some Riemannian metric.

\begin{defin}
\label{def-1.1} 
{\rm
A 
{\it weak chain} 
$\cal V$ 
on
$V$ 
is the following object:

\noindent
a) a sequence of submanifolds of 
$V$:
$V=V_0\supset V_1\supset\ldots\supset V_k\supset V_{k+1}=\emptyset$,
$k\leq {n+1}$, 
${\rm dim}\ V_i = n-i$,
$i=0,\ldots ,k$
(we also set
$V_{-1}=E^{n+N}$);

\noindent
b) some unit vector fields
$v_i$, 
$i=2,\ldots ,k$,
where each
$v_i$
is defined on
$V_i\setminus V_{i+1}$,
is normal to
$V_{i-1}$
in
$V_{i-2}$
and cannot be extended (as such a unit normal vector field) to any subset 
$C\subset V_i$ 
which has a nontrivial  intersection with 
$V_{i+1}$.
}
\end{defin}
\bigskip

We are going to use the notion of a weak chain to 
describe the structure of singularities of tangency of
$V$
to a plane field
${\cal N}$
on
$E$.
Such singularities can be described by an appropriate weak chain equipped
with an additional vector field 
$v_1$
over
$V_1\setminus V_2$.
We are interested in the following two different situations and for each of
these situations
$v_1$
is defined differently.

\smallskip
\noindent
{\sl
Case I:
}
$E=T^\ast L$,
$V$
is a Lagrangian submanifold of 
$E$
and
${\cal N}$
is a 
$\Sigma^2{\hbox{\rm -nonsingular}}$
Lagrangian distribution along
$V$.

\smallskip
\noindent
{\sl 
Case II:
}
$E=M^n\times K^1$,
$V^n\subset E^{n+1}$
is a submanifold,
${\cal N}$
is a line field on
a neighborhood of
$V$
in
$M\times K^1$.
\smallskip

\noindent
We need to consider two cases in order to distinguish later between general
Hamiltonian isotopies of a Lagrangian submanifold and those of them  
which can be covered by a deformation of generating functions with a 
fixed one-di\-men\-si\-onal
space of additional variables (as needed for 
Proposition~\ref{thmC}). 
In the rest of this section we discuss the first case. The second case
will be discussed in Section~\ref{gen-f}. 
Thus in the remainder of this section
$E=T^\ast L$
and
$V$
is a Lagrangian submanifold of
$T^\ast L$.

\begin{defin}
\label{def-1.2}
{\rm
A 
{\it 
chain on 
$V$} 
is a weak chain on
$V$
equipped with an additional unit vector field
$v_1$
defined on the 
{\sl
whole
}
$V_1$
and normal to
$V_1$
in
$V_0$.
}
\end{defin}

\begin{figure}
\centerline{\psfig{figure=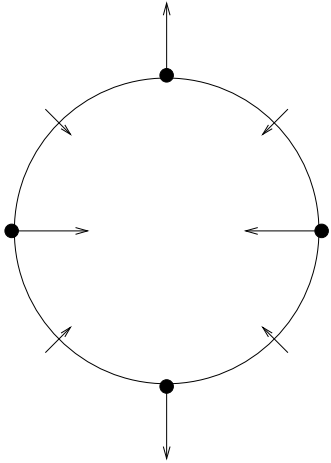,height=70mm}}
\bigskip
\bigskip
\caption{The chain corresponding to the Lagrangian collapse.}
\label{fig1}
\end{figure}

\begin{exam}
\label{exam-1.0}
{\rm
A chain on a 2-dimensional surface is a union of mutually non-in\-ter\-sect\-ing
cooriented curves with some marked points on them. Each marked point is 
equipped with a unit vector normal to the curve this marked point lies on. An example of a
chain in the two-dimensional case is presented on Fig.~\ref{fig1}.
}
\end{exam}

\begin{defin}
\label{def-1.3}
{\rm
Two chains (weak chains) on two submanifolds
$V_1 , V_2 \subset E$
are called
{\it equivalent,}
if there exists an isotopy between
$V_1$
and
$V_2$
in
$E$
which maps one chain (weak chain) into another. 
}
\end{defin}
\bigskip

The following Definition~\ref{def-1.4} 
shows what chains have to do with singularities.
First, we need an auxiliary lemma.

\begin{lem}
\label{lem-1.0}
Let
$x\in V$
be a point.
If
${\cal N}$
is a
$\Sigma^2{\hbox{\it -nonsingular}}$
Lagrangian distribution defined along a neighborhood of
$x$
in
$V$,
then it can be extended to a Lagrangian foliation of a neighborhood of
$x$
in
$T^\ast L$.
\end{lem}

\bigskip
\noindent
{\bf Proof of Lemma~\ref{lem-1.0}.}
If
${\cal N}$
is transversal to
$V$
near
$x$,
then the statement follows from the Lagrangian neighborhood theorem
\cite{We1},\cite{W}.

Otherwise, consider a neighborhood
$U$
of 
$x$
in
$V$
over which
${\cal N}$ 
is defined.
Consider those points of 
$U$
where
${\cal N}$
is not transversal to 
$V$.
The intersection of
${\cal N}$
with 
$T_\ast V$
at each such point is a line. Extend this line field to a 1-dimensional
subbundle 
$l$
of 
${\cal N}$
over
$U$.
The orthogonal complement (with respect to a Riemannian metric) of 
$l$
in
${\cal N}$
is an isotropic bundle 
$l^\perp$
over 
$U$.
Therefore (see 
\cite{We1})
there exist symplectic coordinates 
$({\bf p},{\bf q})$,
${\bf p}= (p_1,\ldots ,p_n)$,
${\bf q}= (q_1,\ldots ,q_n)$,
on a neighborhood of
$x$
in
$T^\ast L$
in which 
$V$
is given by the equation
${\bf q}= 0$
and 
$l^\perp$
(defined over the Lagrangian
${\bf p}{\hbox{\rm -plane}}$)
is given by the equations
${\bf p}=0$, 
$q_n = 0$.
In these new coordinates it is easy to extend the line field
$l$
from the
${\bf p}{\hbox{\rm -plane}}$
to a neighborhood of
$x$
in
$T^\ast L$
so that the plane field spanned by
$l$
and
$l^\perp$
is integrable.
\b
\bigskip

\begin{defin}
\label{def-1.4} 
{\rm
\cite{El}
Let 
${\cal N}$
be a generic  
$\Sigma^2{\hbox{\rm -nonsingular}}$
Lagrangian distribution along a Lagrangian submanifold
$V\subset T^\ast L$.
We define 
{\it
the chain on
$V$
associated with 
${\cal N}$
}
as follows.

First, we define the weak chain. 
For each point of
$V$
we extend 
${\cal N}$
to a Lagrangian foliation of the neighborhood of
$x$
in
$T^\ast L$.
Then using this extension and the fact that
${\cal N}$ 
is generic
we identify 
$V$
and
${\cal N}$
in a neighborhood of
$x$
in
$T^\ast L$
with some local model. Using this local model we determine to which
$V_i$
in the chain the point
$x$
should belong and define the corresponding vector field
$v_i$
at
$x$.
Up to the equivalence of chains, the resulting chain does 
not depend on the extension of
${\cal N}$
and on the identification of a neighborhood of
$x$
in
$T^\ast L$
with the local model.

Let us start with constructing a local model.
Consider a map
$\mu_s:\Rn\to\RnN$,
$s=0,\ldots,n$,
given by the formula
\[\mu_s (t_1,...,t_{n-1},x) =
(t_1,\ldots ,t_{n-1}, \sum_{i=1}^{s-1} t_i x^i + x^{s+1},x,
0,\ldots ,0).\]
Let us denote by
${\cal L}$
an
$N{\hbox{\rm -dimensional}}$
vector bundle on
$\RnN$
formed by the fibers of the projection 
$\pi: \RnN = {\bf R}^n\times {\bf R}^N\to \Rn$.
Let us denote by
$\sigma_s$,
$s=2,\ldots,n$,
the unit vector in
$\RnN$
at
$0\in\RnN$
pointing in the positive direction of
the axis
$t_{s-1}$.

Let us use Lemma~\ref{lem-1.0} 
and extend
${\cal N}$
to a neighborhood of
$x$
in
$T^\ast L$.
Then, according to 
\cite{M}, 
since
${\cal N}$
is generic, for each point
$x\in V$
there exists a neighborhood
$U$
of
$x$
in
$T^\ast L$,
a diffeomorphism
$h:U\to\RnN$
and an integer number 
$s$ 
$(0\leq s\leq n)$
such that
$h(x)=0\in\RnN$,
$h(V\cap U) = \mu_s (\Rn)$
and
$dh$
maps
$\left. {\cal N}\right|_U$
onto
${\cal L}$.

Thus for each point 
$x\in V$
we have a correctly defined number
$s$.
The closure of the set of points for
which the number is equal to
$s$
is a smooth closed
$(n-s){\hbox{\rm -dimensional}}$
submanifold
$V_s$
of
$V$.
Obviously,
$V_0=V$.
We set
$V_{-1}$
to be
$T^\ast L$.

If 
$x\in V_s\setminus V_{s+1}$,
$s\geq 2$,
then the vector
$d(h^{-1})(\sigma_s)$
is tangent to
$V_{s-2}$
and transversal to
$V_{s-1}$.
We define the vector 
$v_s$
at 
$x$
to be a unit vector normal to
$V_{s-1}$
in
$V_{s-2}$
so that it defines the same coorientation of 
$V_{s-1}$
in
$V_{s-2}$
as the vector
$d(h^{-1})(\sigma_s)$.

Thus we have defined the weak chain. 
We define the vector field
$v_1$
to be a unit vector field on
$V_1$
normal to
$V_1$
in
$V$
and coorienting
$V_1$
in
$V$
according to the canonical Maslov coorientation 
of 
$V_1$
in
$V$
(see \cite{Ar4}).
These conditions uniquely define the vector field
$v_1$.

Thus we have defined a weak chain and a vector field
$v_1$.
One can check that it is, indeed, a chain. Up to the 
equivalence of chains this chain does not depend on the choices made in
the construction.
}
\end{defin}

\begin{rem}
\label{kart-sborka}
{\rm
In the two-dimensional case the way in which we associate a vector
$v_2$
to each cusp (to each point in
$V_2$)
can be seen on Fig.~\ref{fig2}.
}
\end{rem}

\begin{rem}
\label{rem-1.1}
{\rm
Consider the Lagrangian distribution
$Vert$
along a Lag\-ran\-gian submanifold
$V\subset T^\ast L$
which associates to each point
$x\in V$
the tangent space at
$p$
to the fiber of the projection
$\pi: T^\ast L\to L$
which contains
$x$. 
Assume that the Lagrangian distribution
$Vert$
is generic and
$\Sigma^2{\hbox{\rm -nonsingular}}$. 
Consider a weak chain on
$V$
associated with
$Vert$
as in 
Definition~\ref{def-1.4}.
In this weak chain the set
$V_i\setminus V_{i+1}$
consists of those points in
$V$
where the map
$\left.\pi\right|_V: V\to L$
has a singularity of the Thom-Boardman type 
${\Sigma}^{\underbrace{11\ldots 1}_{i\ \rm times}}$.

In fact, this weak chain is completely determined by the map
$\pi: V\to L$
(one can actually associate a weak chain to any smooth map
$f:V^n\to L^n$
which has only 
${\Sigma}^{11\ldots 1}$
singularities 
\cite{El}).
}
\end{rem}
\bigskip

When a Lagrangian distribution
${\cal N}$
along
$V$
changes, some bifurcations may happen to the chain on
$V$
associated with
${\cal N}$.
There are two sorts of such typical bifurcations:
{\it 
direct surgeries
}
and
{\it
inverse surgeries
}.
Let us start with defining a direct surgery as a formal operation on
a chain.
We first define a necessary setup needed to perform a direct surgery.

\begin{defin}
\label{def-1.7} 
{\rm
\cite{El}
A
{\it 
basis of a direct surgery}
of order
$s\geq 1$ 
and index
$p\geq 0$
is formed by a chain
$\cv$
and an embedding
$\phi :D^p\times [-1,1]\to V_{s-1}$.
When 
$p=0, s=n$
this embedding must be equipped with a
pair of nonzero vectors
$\nu_1 ,\nu_2$
which are normal to
$V_{n-1}$
in
$V_{n-2}$
and which have opposite directions; 
we also assume in this case that
$\phi (D^0\times [-1,1])\subset V_{n-1}\setminus V_{n-2}$.
When 
$p=0, s=1$,
then, in addition to the data above, the basis includes
a pair of nonzero vectors
$\nu_3 , \nu_4$
tangent to the curve
$\phi :D^0\times [-1,1]\to V_1$
at the two endpoints and pointing either both inside the curve or both 
outside the curve.

All this data has to satisfy the following conditions:

\smallskip
\noindent
a) 
$\phi (D^p\times 0)\cap V_s=
\phi (D^p\times [-1,1])\cap V_{s+1}=
\phi (\partial D^p\times 0)$
and the intersections are normal;

\noindent
b)
$\phi (D^p\times [-1,1])\cap V_s=\phi (\partial D^p\times [-1,1])$
and the intersection is normal;

\noindent
c) the vector field 
$\left. v_{s+1}\right|_{\phi (\partial D^p\times 0)}$
points inside the disk
$\phi (D^p\times 0)$;

\noindent
d) If 
$s>1$
the vector fields
$\left. v_s\right|_{\phi (\partial D^p\times 1)}$
and
$\left. v_s\right|_{\phi (\partial D^p\times \{-1\})}$
can be extended over
$\phi (D^p\times 1)$
and, respectively over
$\phi (D^p\times \{-1\})$,
as vector fields normal to
$V_{s-1}$
in
$V_{s-2}$.

\noindent
e) 
If
$s=1$, 
the vector field
$\left. v_1\right|_{\phi (\partial D^p\times 0)}$
points either everywhere inside  or everywhere outside the connected
component of
$V_0\setminus V_1$
which contains
$\phi (D^p\times 0)$.
}
\end{defin}
\bigskip

We recall the definition of a 
{\it
direct surgery.
}

\begin{defin}
\label{def-1.8}
{\rm
\cite{El}
Suppose that a chain
$\cv$ 
and an embedding
$\phi :D^p\times [-1,1]\to V_{s-1}$
(possibly with a pair of vectors
$\nu_1 ,\nu_2$
or
$\nu_3 ,\nu_4$) 
form a basis of a direct surgery of order
$s$
and index 
$p$.
The result of the direct surgery associated with the basis
is a new chain
${\cal V}^\prime$
defined as follows.

\smallskip
For
$i\neq s,s+1$
let us put
$V_i^\prime = V_i$
and let us define the pair 
$V_s^\prime \supset V_{s+1}^\prime$
of submanifolds of the manifold
$V_{s-1}^\prime = V_{s-1}$
as the result of the simultaneous attachment of a handle 
of index 
$p$
to the pair
$V_s\supset V_{s+1}$
along
$\phi (D^p\times 0)$
in
$V_{s-1}$.
In other words, we do the following. 
We consider a (closed) tubular neighborhood 
$T$
of
$\phi (D^p\times 0)$
in  
$V_{s-1}$
(more precisely, in a normal bundle of
$\phi (D^p\times 0)$
in  
$V_{s-1}$)
and present it as a product
$T=T_1\times [-1,1],$
where

\noindent
(i) 
$T_1$
is in general position with respect to
$\phi (D^p\times [-1,1])$
in
$V_{s-1} ;$ 

\noindent
(ii)
$T_1\cap V_s =T\cap V_{s+1}.$

We put 
$N = T\cap V_s$,
$N_1 = T\cap V_{s+1} = T_1\cap V_s$
and define
$V_s^\prime$
as
$V_s^\prime = (V_s\setminus N)\cup \partial T$.
If 
$N_1\neq\emptyset$ 
and
$(V_{s+1}\setminus N_1)\cap \partial T_1=\emptyset$,
then we put
$V_{s+1}^\prime = V_{s+1}\setminus N_1$,
else we put 
$V_{s+1}^\prime = (V_{s+1}\setminus N_1)\cup \partial T_1$.
Finally we smoothen the corners to get smooth manifolds
$V_s^\prime$
and
$V_{s+1}^\prime$.
\smallskip

The vector fields 
$v_i^\prime$
are defined as follows.
For 
$i\neq s-1,s,s+1$
we set
$v_i^\prime = v_i$.
We can define
$v_{s-1}^\prime$
(and in the case when
$p>0$
also
$v_s^\prime ,v_{s+1}^\prime$)
uniquely (up to the equivalence
of chains), if we require that 
$v_{s-1}^\prime ,v_s^\prime ,v_{s+1}^\prime$
coincide with 
$v_{s-1} ,v_s , v_{s+1}$
respectively outside of a neighborhood of
$\phi (D^p\times [-1,1])$
in
$V$
and that together with the already constructed objects
$V_i^\prime$ 
and
$v_i^\prime$
they form a chain.

If 
$p=0$,
then the surgery adds to
$V_s$
another connected component
$\partial T$
and to
$V_{s+1}$
another connected component
$\widetilde V\subset \partial T$.
When 
$1<s<n$
we define
$v_s^\prime , v_{s+1}^\prime$
by requiring in addition that the field
$\left. v_{s+1}^\prime\right|_{\widetilde V}$
points outside
$T$.
This additional requirement together with the requirements above
determines 
$v_s^\prime , v_{s+1}^\prime$
uniquely (up to the equivalence of chains).
If
$s=n$,
then
$\widetilde V=\emptyset$
and we define
$v_s^\prime$
at the points
$\phi (D^0\times \{-1\})$
and
$\phi (D^0\times 1)$ 
as
$\nu_1$
and
$\nu_2$
respectively. 
Finally, if
$p=0 , s=1$
we define
$v_1$
on
$\partial T$
as a vector field normal to
$\partial T$
in
$V_0$
which agrees with
$\nu_3 , \nu_4$
at the endpoints of
$\phi (D^0\times [-1,1])$

The chain we obtain is unique up to the equivalence of chains.
}
\end{defin}
\bigskip

From now on, whenever we speak about a direct surgery along
a disk
$D^p$,
we omit
$\phi$
and identify
$\phi (D^p)$
with
$D^p$
and 
$\phi (D^p\times[-1,1])$
with
$D^p\times[-1,1]$.

An 
{\it 
inverse surgery,
} 
as a formal operation on a chain,
can be defined in a spirit similar to the definition of a direct surgery
(see \cite{El} for details).
Inverse surgeries are exactly the operations on chains 
inverse to direct surgeries.

The following proposition is analogous to a result mentioned in 
\cite{El}
and can be proved (as the result in 
\cite{El}) 
by the methods developed in
\cite{Ar3} 
(also see
Section 22.4 
of 
\cite{AGV}).

\begin{prop}
\label{prop-1.1} 
Let
$\{ {\cal N}^t\}$
be a generic 1-parameter family of
\break
$\Sigma^2{\hbox{\it -nonsingular}}$
Lagrangian distributions along
$V$ 
(${\rm dim}\  V\geq 2$). 
Then the chain on
$V$
associated with
${\cal N}^1$
can be obtained from the chain
associated with
${\cal N}^0$
by a sequence of direct and inverse surgeries.
\end{prop}

\begin{figure}
\centerline{\psfig{figure=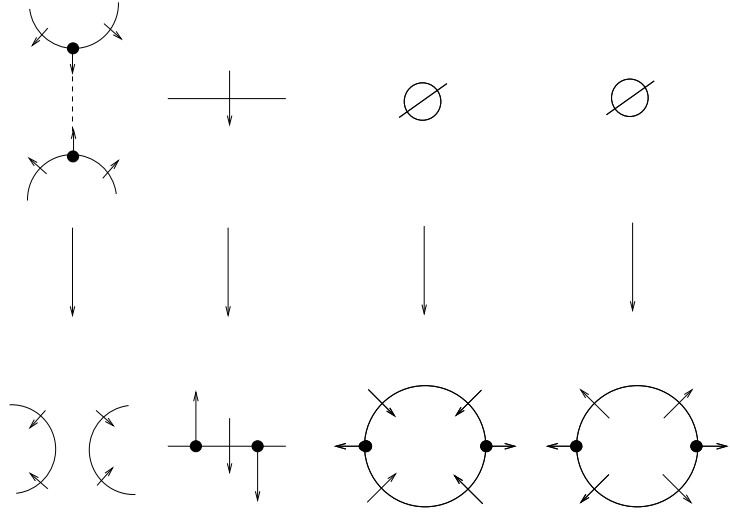,height=80mm}}
\caption{The list of possible direct surgeries for the two-dimensional case.}
\label{fig3a}
\end{figure}

\begin{figure}
\centerline{\psfig{figure=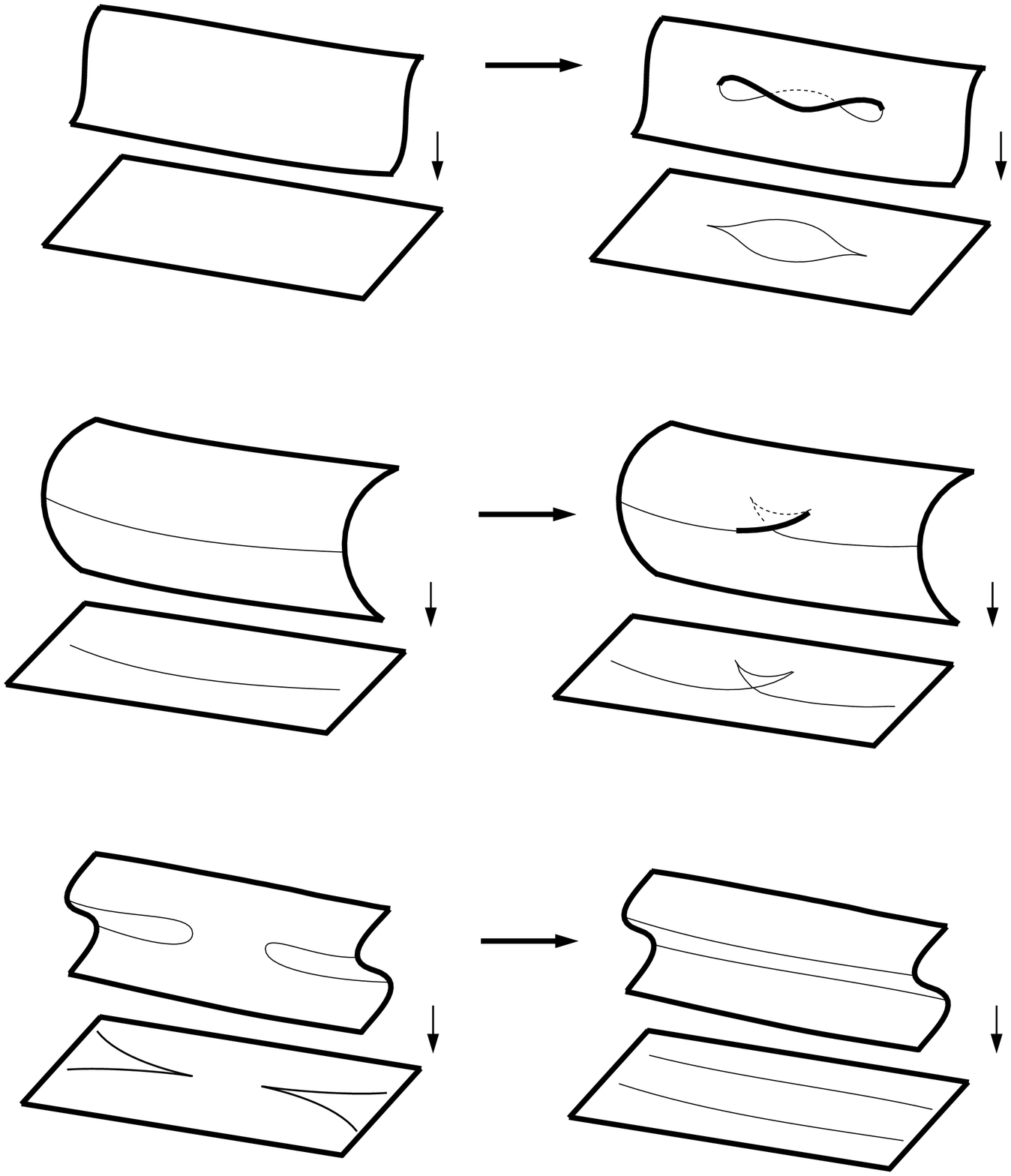,height=175mm}}
\caption{Direct surgeries viewed as changes of surfaces with respect to a 
projection.}
\label{fig3b}
\end{figure}

\begin{exam}
\label{exam-1.1}
{\rm
All possible direct surgeries for the two-dimensional case are presented 
on Fig.~\ref{fig3a}. The actual changes of a surface (with respect to a fixed
vertical line field) are presented on Fig.~\ref{fig3b}.
The possible inverse surgeries are exactly the opposite operations.
}
\end{exam}

\begin{prop}
\label{prop-1.2}
{\rm
\cite{El}
}
Consider an inverse surgery on a chain 
\break
${\cal V} = \{ V_i,v_i\}$
which transforms 
${\cal V}$
into another chain
${\cal V}^\prime = \{ V_i^\prime, v_i^\prime\}$.
Suppose that the inverse surgery is an inverse of a direct surgery 
of order
$s$
along a disk
$D^p$
(i.e. this direct surgery transforms
${\cal V}^\prime$
into
$\cal V$).
Suppose also that there exists an embedding 
$\alpha : D^1\to V_0$
which satisfies the following properties:

\noindent
(i)
$\alpha (0)\in \partial D^p$,
$\alpha (1)\not\in \partial D^p$,
$\alpha (D^1)\cap V_1 = \alpha (\partial D^1)$
and 
$\alpha (D^1)$ 
is normal to
$V_1$
at 
$\alpha (\partial D^1)$;

\noindent
(ii)
the vectors from the vector field
$v_1$
at the endpoints of
$\alpha (D^1)$
either both point inside 
$\alpha (D^1)$
or both point outside of it.

\smallskip
\noindent
Then the inverse surgery can be represented as a
composition of direct surgeries.
\end{prop}

\begin{defin}
\label{def-1.9}
{\rm
\cite{El}
We call an inverse surgery satisfying the hypothesis of 
Proposition~\ref{prop-1.2}
{\it 
complete.
}
}
\end{defin}

\begin{figure}
\centerline{\psfig{figure=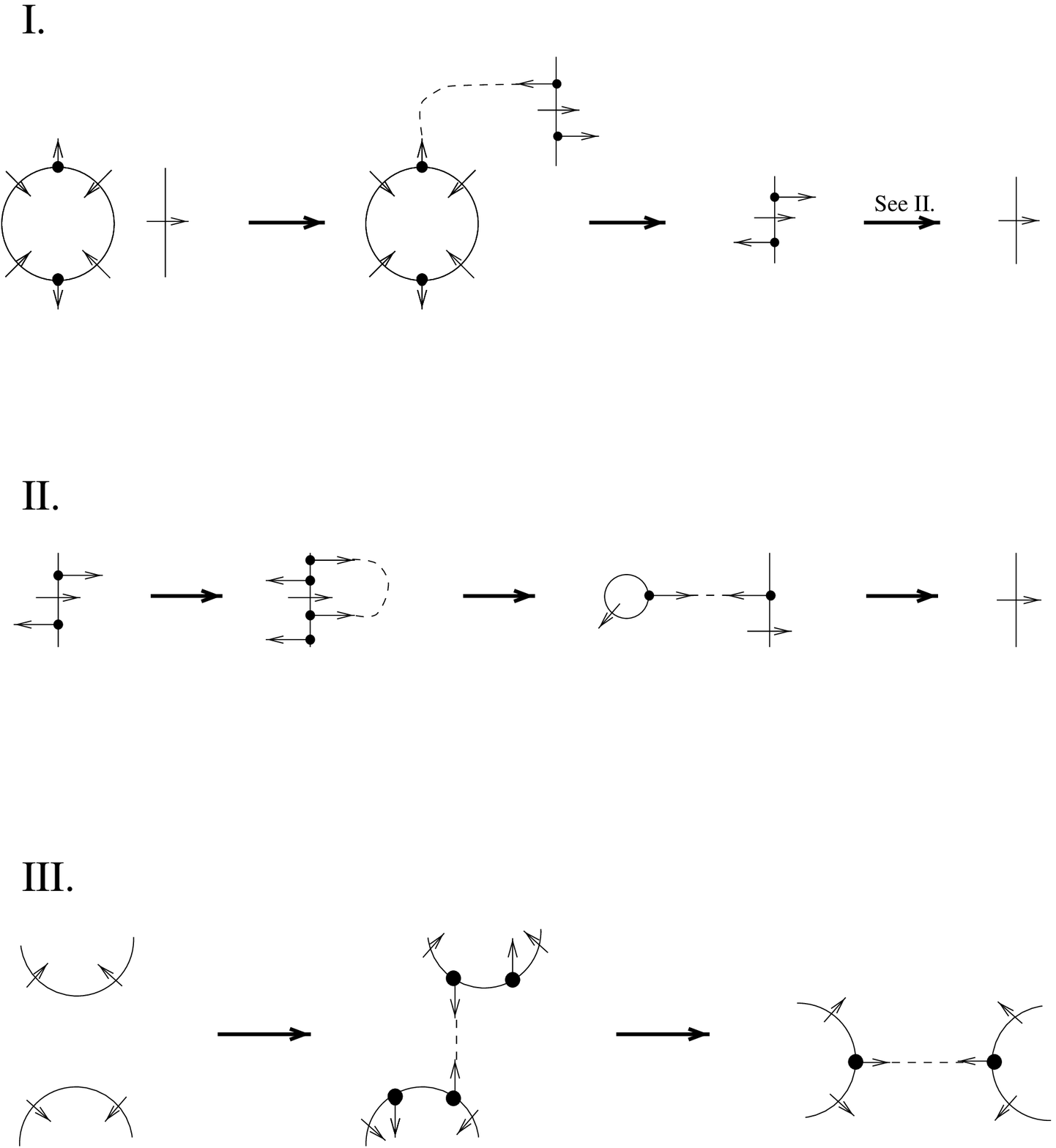,height=175mm}}
\caption{Complete inverse surgeries presented as compositions of direct surgeries 
in the two-dimensional case.}
\label{fig4}
\end{figure}

\begin{exam}
\label{exam-1.99}
\rm{
The proof of Proposition~\ref{prop-1.2} 
in the two-dimensional case is presented on 
Fig.~\ref{fig4}.
}
\end{exam}

\section{Proof of Theorem B}
\label{sect-pf-thmB}

In this section we state that a direct surgery can be realized
by a Hamiltonian isotopy. We also state that by a Hamiltonian isotopy
we can create additional double fold pairs of  
spheres in a neighborhood of a point in
$L$
at which the foliation
${\cal F}$
is transversal to
$L$.
The proofs of the statements are contained in the subsequent sections.
From these
results we deduce Theorem B.

\bigskip
\noindent
{\bf Notational convention.}

\smallskip
\noindent
From now on we do not distinguish between
$L$
and the zero-section of
$T^\ast L$.
We denote by the same letter
${\cal F}$
the foliation on a neighborhood 
${\cal U}$
of the zero-section
$L$
in
$T^\ast L$
and the Lagrangian distribution along
$L$
induced by this foliation.
\bigskip

\begin{prop}
\label{prop-2.1}
Let 
$\cal V$ 
be the chain on 
$L$
associated with
${\cal F}$ 
and let 
${\cal V}_1$ 
be the result of a direct surgery of order 
$s$,
$s\geq 1$, 
along a disk 
$D^p\subset V_{s-1}$. 
Then for any neighborhood 
$U$ 
of 
$D^p$ 
there exists a Hamiltonian flow 
$\{ h_t\}$,
$0\leq t\leq 1$, 
such that:

\smallskip
\noindent
(i) 
the chain 
${\cal V}_1$ 
is equivalent to the chain on the Lagrangian submanifold
$L_1 = h_1(L)$
associated with
${\cal F}$;

\noindent
(ii)  
$h_t = id$
outside 
$U$
for any
$t$,
$0\leq t\leq 1$.

\smallskip
\noindent
The Hamiltonian isotopy 
$\{ h_t\}$
above can be chosen in such a way that it is 
compactly supported and arbitrarily
$C^0{\hbox{\it -small}}$
and so that 
$\displaystyle \max_{t,\gamma}
\mid \int_{h_t(\gamma)} pdq \mid$ 
is also arbitrarily
$C^0{\hbox{\it -small}}$.

\end{prop}
\bigskip

We will prove 
Proposition~\ref{prop-2.1} in 
Section~\ref{sect-gen-f-pfs}.

\begin{prop}
\label{prop-2.2}
Let
$B\subset L$
be a ball in
$L$
such that
${\cal F}$
is transversal to
$L$
along 
$B$.
Let
${\cal V}$
be the chain on
$L$
associated with
${\cal F}$.
Consider a new chain
${\cal V}_1$
on
$L$
which is a union of
${\cal V}$
and a double fold pair of spheres lying in
$B$.

Then there exists a Hamiltonian isotopy 
$\{ h_t (L)\}$,
$0\leq t\leq 1$,
of
$L$
such that the chain
${\cal V}_1$
is equivalent to the chain on
$L$
associated with the Lagrangian distribution
$(h_1)^{-1}_\ast ({\cal F})$
along
$L$.

The Hamiltonian isotopy
$\{h_t\}$
can be chosen in such a way that it is 
compactly supported and arbitrarily
$C^0{\hbox{\it -small}}$
and so that 
$\displaystyle \max_{t,\gamma}
\mid \int_{h_t(\gamma)} pdq \mid$ 
is also arbitrarily
$C^0{\hbox{\it -small}}$.

Any of the two possible combinations of Maslov coorientations of 
$S_1$
and
$S_2$
can be realized by a Hamiltonian isotopy as above. 
\end{prop}
\bigskip

We will prove Proposition~\ref{prop-2.2} in Section~\ref{double-folds-pfs}.

\bigskip
\noindent
{\bf Proof of Theorem B.}
Let 
$\{ {\cal N}^t\}$
be a deformation connecting
${\cal F}$
with
${\cal N}$
in the class of 
$\Sigma^2{\hbox{\rm -nonsingular}}$
Lagrangian distributions along
$L$.
We apply Proposition~\ref{prop-1.1} 
and get a sequence of direct and inverse surgeries associated to the homotopy
$\{ {\cal N}^t\}$. 
Then we deform
$L$
step by step by  
$C^0{\hbox{\rm -small}}$
Hamiltonian isotopies so that the chain on the deformed submanifold
(which we identify with 
$L$)
associated with
${\cal F}$
changes according to the given surgery.

Namely, let us take the first surgery from the sequence.
If it is 
{\sl
a direct surgery,
} 
we realize it by a
$C^0{\hbox{\rm -small}}$
Hamiltonian isotopy with the help of
Proposition~\ref{prop-2.1}. 

If it is
{\sl
a complete inverse surgery, 
}
then we use 
Proposition~\ref{prop-1.2} 
and decompose this inverse surgery into a 
product of direct surgeries reducing the problem
to the previous situation.

If this is an inverse surgery which is
{\sl
not complete, 
}
we make 
it complete by the following trick. We can pick a point
$x\in L$ 
at which
${\cal F}$
is transversal to
$L$
and, using 
Proposition~\ref{prop-2.2}, 
create (by a 
$C^0{\hbox{\rm -small}}$ 
Hamiltonian isotopy supported in an arbitrarily
small neighborhood of
$x$)
a double fold along some spheres
$S_1$
and
$S_2$.
We can assume then that the result of the inverse surgery is the chain  
${\cal V}_1$ 
formed by
${\cal V}$ 
and the spheres
$S_1$
and
$S_2$.
By an appropriate choice of
$x$
and of the Maslov coorientations orientations of 
$S_1$
and
$S_2$
(see Proposition~\ref{prop-2.2})
we can always make the inverse surgery complete and reduce the problem to
the previous case.

Thus we can deform
$L$
by a 
$C^0{\hbox{\rm -small}}$
Hamiltonian isotopy into a new Lagrangian submanifold
$L^\prime\subset {\cal U}\subset T^\ast L$. 
Consider a small tubular neighborhood of
$L^\prime$
in
${\cal U}$. 
This tubular neighborhood is  
symplectomorphic to a neighborhood
of the zero-section in
$T^\ast L$
by the Lagrangian neighborhood theorem 
\cite{We1},\cite{W}.
Therefore we can identify
$L^\prime$
with
$L$,
take the next surgery from the sequence
and repeat the whole procedure for the new 
$L$.
Then we take the next surgery from the sequence etc.

Thus working step by step, we construct a sequence of compactly supported 
and arbitrarily
$C^0{\hbox{\rm -small}}$
Hamiltonian isotopies of
$L$.
According to Proposition~\ref{prop-2.2},
each of the Hamiltonian isotopies can be chosen so that, if we fix a path 
$\gamma$
in
$L$
and consider its deformation 
$\{\gamma_t\}$
under the isotopy, then the number
$\mid\int_{\gamma_t} pdq \mid$
is bounded by a uniform arbitrarily small number which does not depend on
$t$
and
$\gamma$.
The composition of all the isotopies above gives us the necessary
Hamiltonian isotopy
$\{ h_t(L)\}$. 
The chain on
$h_1 (L)$
associated with
${\cal F}$
is equivalent to the union of the chain
${\cal V}$
associated with 
${\cal N}$
and possibly a number of additional double fold pairs of spheres 
that we might have
created when we decomposed an incomplete inverse surgery into a product
of direct surgeries.
\b
\bigskip

\section{The two-dimensional case}
\label{two-dim-case}

\subsection {Preliminaries}
\label{prelim}

Let us start with a discussion of the Grassmannian 
$\Lambda_2^{+}$
of oriented Lagrangian planes in the standard symplectic space
${\bf R}^4$.
This Grassmannian can be represented as the quotient
$U(2)/SO(2)$
which is diffeomorphic to
$S^2\times S^1$
(see e.g.
\cite{Ar4} 
or
\cite{AG}).

Let us fix a Lagrangian plane
$o$
in
${\bf R}^4$.
Let
$o^{+}$,
$o^{-}$
be the two points in
$\Lambda^{+}_2$
corresponding to the plane 
$o$
taken with the two possible orientations.
The points in
$\Lambda^{+}_2$
corresponding to oriented Lagrangian planes
which are not transversal to
$o$
in
${\bf R}^4$
form a 
$2{\hbox{\rm -dimensional}}$
surface
$\Sigma\subset \Lambda^{+}_2$
which is smooth outside 
$o^{+}$
and
$o^{-}$
and which has  isolated singularities at
$o^{+}$
and
$o^{-}$
topologically equivalent to the Morse singularity of the cone
$\{(x,y,z)\in {\bf R}^3\  |\   x^2+y^2-z^2=0\}$
at zero. 
Topologically the surface
$\Sigma$
is formed by some two spheres
$S_{+}$
and
$S_{-}$
glued to each other at the two points
$o^{+}$
and
$o^{-}$.
The homology classes of
$S_1$
and 
$S_2$
are equal (up to sign) to the homology class of
$S^2\times \{pt\}$
in
$\Lambda^{+}_2\cong S^2\times S^1$.
The complement
$\Lambda^{+}_2\setminus\Sigma$
consists of two connected components homeomorphic to a 
three-dimensional ball. These two components, the 
``positive'' one and the ``negative'' one, correspond to 
the two possible orientations of
${\bf R}^4$:
a point in the ``positive'' (resp. the ``negative'') component 
corresponds to an oriented  Lagrangian plane transversal to
$o$
which together with
$o^{+}$
defines the positive (resp. the negative) orientation of
${\bf R}^4$
with respect to the symplectic form
$\omega$.
The canonical Maslov coorientation of
$\Sigma\setminus \{o^{+},o^{-}\}$
(see \cite{Ar4})
points from the ``positive'' component to the ``negative'' on one of the
spheres
$S_{+}$,
$S_{-}$
and from the ``negative'' component to the ``positive'' component on another.
Locally,
$\Sigma$
subdivides a small neighborhood of
$o^{-}$
(resp.
$o^{+}$)
into three connected components. One of them (a ``big'' one) has a common 
boundary with the two others and contains the oriented Lagrangian planes which 
together with 
$o^{-}$ 
(resp.
$o^{+}$)
determine the positive orientation of
${\bf R}^4$.

Another way to describe
$\Sigma$
is the following. Take a line
$l$
in the plane
$o$.
The oriented Lagrangian planes in
${\bf R}^4$
containing
$l$
form a simple closed curve
$\gamma_l\subset\Sigma$
in the Grassmanian
$\Lambda_2^{+}$.
For any 
$l\subset o$
the curve
$\gamma_l$
passes through
$o^{+}$
and
$o^{-}$.
The surface
$\Sigma$
is formed by a revolving curve
$\gamma_l$
as the line
$l$
revolves (counterclockwise) in
$o^{+}$,
and a full turn of
$l$
in
${\bf RP}^1$
corresponds to a full turn of
$\gamma_l$
in 
$\Sigma$.
By specifying that we rotate
$l$
counterclockwise in
$o^{+}$,
we have chosen some generators
$g_1, g_2$
of
$\pi_1 (S_1\setminus\{o^{+},o^{-}\})$,
$\pi_1 (S_2\setminus\{o^{+},o^{-}\})$.
This also determines a specific orientation for each curve
$\gamma_l$.

Each curve
$\gamma_l$
from above represents a generator  
$C_1$
of the fundamental group
$\pi_1(\Lambda_2^{+})\cong {\bf Z}$.
Moreover, the curve 
$\gamma_l$
is oriented and so we can fix a specific generator
$c_1\in H^1(\Lambda_2^{+})$
as a cohomology class dual to
$C_1$.
We define an orientation of
$\Lambda_2^{+}$
by requiring that the intersection number of 
$C_1$
with the oriented disk bounded by 
$g_i$
is
$1$.
The disk bounded by
$g_i$
is oriented in such a way that
$g_i$
rotates counterclockwise with respect to the orientation of the disk.
Another representative of 
$C_1$
can be obtained by taking any plane
$o_1\in \Lambda_2^{+}$
and taking the closed path
$\varsigma_{o_1} (t) = e^{\pi i t} o_1$,
$0\leq t\leq 1$.

Let 
$c_2\in H^2(\Lambda_2^{+})$
be the cohomology class dual to the homology class
$S^2\times \{pt\}$
in
$S^2\times S^1$,
where the orientation of 
$S^2\times \{pt\}$
is chosen so that the intersection number of
$S^2\times \{pt\}$
with 
$\{pt\}\times S^1$
(oriented as above)
is equal to
$1$
with respect to the orientation of the whole space
$\Lambda_2^{+}$
chosen as above.

Consider a (constant) oriented Lagrangian plane field
$o^{+}$
on
${\bf R}^4$.
Let
$L^2\subset {\bf R}^4$
be an (immersed) oriented Lagrangian surface with boundary such that
$L$
is transversal to
$o$
over
$\partial L$.
Let 
$f: L\to \Lambda_2^{+}$
be the Gaussian map of
$L$.
Then
$c=f^\ast c_1\in  H^1(L,\partial L)$
is the Maslov class of 
$L$
with respect to
$o$
(or with respect to the Lagrangian projection of
${\bf R}^4$
along
$o$)
-- see \cite{Ar4}.
Also observe that since 
$L$
is oriented each component of
$\partial L$
acquires a canonical orientation and therefore we 
have a canonical trivialization of
$T_\ast L$
near 
$\partial L$.
This allows us to define the (relative) Euler class 
$e\in H^2 (L,\partial L)$
of the oriented plane bundle
$o^{+}$
on
$L$.
One can check that with our choice of orientations the class
$e$
is induced by 
$f$
from the generator
$c_2$   
(see \cite{Bor} or \cite{Vas}, p.239).
More precisely, the trivialization
of
$T_\ast L$
along the boundary
$\partial L$
correctly defines a map
$f^\ast: H^2(\Lambda^{+})\to H^2(L/{\partial L})\cong H^2(L,{\partial L})$.
This finishes the discussion of the Grassmannian
$\Lambda_2^{+}$.

Now we want to extend the definitions of Maslov and Euler classes to the case
of an oriented Lagrangian distribution
${\cal N}$
along an oriented surface (with boundary)
$L$
identified with the zero-section of
$T^\ast L$.
{\sl
As before we consider only those Lagrangian distributions which are
transversal to
$L$
near
$\partial L$.
}
In particular, this implies that there exists a
unique  almost complex structure
$J$
defined near
$\partial L$
on
$T^\ast L$
so that (near
$\partial L$)
one has
$J (T_\ast L) = {\cal N}$
and so that 
$J$
is tamed by the canonical symplectic form
$\omega = d(pdq)$
on
$T^\ast L$.
(We say that an almost complex structure
$J$
on
$T^\ast L$
is 
{\sl
tame}, 
or is 
{\sl
tamed by the symplectic form 
}
$\omega$,
if
$\omega (Jv,v) > 0$ 
and
$\omega (Jv, Jv) = \omega (v,v)$
for any nonzero
$v\in T_\ast (T^\ast L)$).
This almost complex structure 
$J$,
defined near 
$\partial L$,
can be extended to a tame almost complex structure  on the whole
$T^\ast L$.
Together 
$J$
and
$\omega$
define a Riemannian metric
$g$,
where
$g(x,y) = \omega (x, Jy)$,
$x,y\in T_\ast (T^\ast L)$.

If we consider the canonical trivialization of
$T_\ast L$
near
$\partial L$,
we can use
$J$
to find the canonical trivialization of
${\cal N}$
near
$\partial L$.
Therefore there exists a correctly defined (relative) Euler class 
$e({\cal N})\in H^2 (L,\partial L)\cong {\bf Z}$
of
$\cal N$.

To define the Maslov class of 
${\cal N}$
in this situation one should imitate the definition for the case of 
a Lagrangian submanifold of
${\bf R}^{2n}$ 
in the following way.
Take a point 
$x\in L$
and choose some 
$g{\hbox{\rm -orthogonal bases}}$ 
in the tangent space 
$T_x L$
to the zero-section
$L\subset T^\ast L$
and in the Lagrangian plane
${\cal N}_x$.
There exists a unique unitary linear operator 
$u(x)$
on the complex Hermitian space
$T_x (T^\ast L)$
which maps the basis of
$T_x L$
into the basis of
${\cal N}_x$.
(This operator 
$u(x)$
is an element of the fiber over
$x$
of the principal
$U(2){\hbox{\rm -bundle}}$
over
$T^\ast L$
defined by
$J$
and
$g$).
The map
$f(x) = {\rm det}^2\ (u(x))$
from
$L$
to the circle
${\bf S}^1={e^{i\phi}}\subset {\bf C}$
is correctly defined and does not depend on the choice of the 
bases above.
Observe that since 
${\cal N} = J (T_\ast L)$,
the map 
$f$
is identically equal to 1 near 
$\partial L$.
Let
$c^\prime$
be the generator of
$H^1 ({\bf S}^1, pt)$
corresponding to the counterclockwise orientation of
${\bf C}$.
The cohomology class
$c({\cal N}) = f^\ast (c^\prime)\in H^1 (L,\partial L)$
is called the 
{\it
Maslov class
}
of the Lagrangian distribution
${\cal N}$
along
$L$.

The Maslov class of
${\cal N}$
can also be seen in the following way.
As we mentioned before, a generic Lagrangian distribution along
$L$
is
$\Sigma^2{\hbox{\rm -nonsingular}}$
because  Lagrangian singularities of the Thom-Boardman type
$\Sigma^2$
generically have codimension 3
(see \cite{Ar1},\cite{AGV}).
Assume now that
$\cal N$ 
is generic and consider the chain on
$L$
associated with 
${\cal N}$
(see Definition~\ref{def-1.4}).
The submanifold
$V_1$
from this chain, cooriented by the vector field
$v_1$,
defines a cohomological class
$c({\cal N})\in H^1(L,\partial L)$,
which is nothing else but the
Maslov class of 
$\cal N$.

\subsection{Auxiliary results}
\label{auxil}

Let
$L$
be an oriented two-dimensional compact surface (possibly with boundary).

\begin{prop}
\label{prop-22.1}
Let
${\cal N}^1$
and
${\cal N}^2$
be two oriented Lagrangian distributions coinciding near
$\partial L$.
Then
${\cal N}^1$
and
${\cal N}^2$
are homotopic (rel.
$\partial L$)
as oriented Lagrangian distributions if and only if their Euler
and Maslov classes are equal.
\end{prop}

Let
$\cal N$
be a generic
$\Sigma^2{\hbox{\rm -nonsingular}}$
oriented Lagrangian distribution along
$L$.
Consider the submanifold
$V_1$
from the chain on
$L$
associated with
${\cal N}$.
Since
$L$
and
${\cal N}$
are oriented,
$V_1$
subdivides
$L$
into two parts
$L_1$ 
and
$L_2$
so that over
$L_1$
(resp. 
$L_2$)
the oriented submanifold
$L$
and
${\cal N}$
together define a positive 
(resp. the negative)
orientation of
$T^\ast L$.

\begin{prop}
\label{prop-22.2}
{\rm 
\cite{Ha},\cite{El2} 
}
In the situation as above
let us assume that there are
$n_1$
points in 
$V_2$
at which the vectors
$v_2$
point inside
$L_1$
and
$n_2$
points where
$v_2$
points inside
$L_2$.
Then
$\chi (L_1) - \chi (L_2) + n_1 - n_2 = e({\cal N})$.
\end{prop}

\begin{prop}
\label{prop-22.3}
A chain 
$\cal V$
on
$L$
can be realized as a chain associated with some
$\Sigma^2{\hbox{\it -nonsingular}}$
oriented Lagrangian distribution along
$L$
if and only if it satisfies the following property:

\smallskip
\noindent
the submanifold
$V_1$
from the chain 
$\cal V$
subdivides
$L$
into two different (not necessarily connected) parts with the common
boundary
$V_1$.

\end{prop}

\bigskip
\noindent
{\bf Proof of 
Proposition~\ref{prop-22.1}.}
Suppose that
${\cal N}^1$
and
${\cal N}^2$
are homotopic rel. 
$\partial L$.
It follows from the definition that the Maslov classes of
${\cal N}^1$
and
${\cal N}^2$
are equal.
The Euler classes also have to be equal since 
${\cal N}^1$
and
${\cal N}^2$
are isomorphic rel.
$\partial L$.

Now suppose that
$c({\cal N}^1)= c({\cal N}^2)$,
$e({\cal N}^1)= e({\cal N}^2)$.
Consider a vector bundle 
$E\to L$
which associates to each point 
$x\in L$
the tangent space
$T_x (T^\ast L)$.
Consider a fiber bundle
$\Lambda_2^{+} (E)\to L$ 
which associates to each point 
$x\in L$
the Grassmannian of oriented Lagrangian planes in the fiber of
$E$
over
$x$.
Any oriented Lagrangian distribution 
${\cal N}$
on
$L$
is a section of the fiber bundle
$\Lambda_2^{+} (E)\to L$.
We triangulate
$L$
into small simplices so that over each 2-simplex the vector bundle
$E\to L$
is trivial. It follows from Section~\ref{prelim} and from 
the standard obstruction theory that the equality
$c({\cal N}^1)= c({\cal N}^2)$
allows to extend a homotopy between
${\cal N}^1$
and
${\cal N}^2$
from the 0-skeleton of the triangulation to the 1-skeleton 
(keeping it identical near
$\partial L$),
and the equality
$e({\cal N}^1)= e({\cal N}^2)$
allows to extend a homotopy  between
${\cal N}^1$
and
${\cal N}^2$
from the 1-skeleton to the 2-skeleton
(again keeping it identical near
$\partial L$).
\b
\bigskip

\bigskip
\noindent
{\bf Proof of 
Proposition~\ref{prop-22.2}.}
Take some small tubular neighborhoods in
$L$
of each closed curve 
$\gamma_i\subset V_1$
and cut 
$L$
open along the boundaries of those tubular neighborhoods. Over each piece
${L^\prime_i}$
that does not intersect
$V_1$
the plane bundle
${\cal N}$
is isomorphic to
$\pm T_\ast {L^\prime_i}$
(the sign depends on the orientations, i.e. whether 
${L^\prime_i}\subset L_1$
or
${L^\prime_i}\subset L_2$).
Therefore,
$e(\left. {\cal N}\right|_{L^\prime_i}) =\pm\chi ({L^\prime_i})$.
Orient each
$\gamma_i$
in such a way that, if we rotate a positive tangent vector to
$\gamma_i$ 
by 90 degrees counterclockwise, we get a normal vector
which points into
$L_1$.
It follows from the discussion in Section~\ref{prelim}
that over each cylindrical piece
$L_i^{\prime\prime}$
containing a curve
$\gamma_i\subset V_1$
the Euler class
$e(\left. {\cal N}\right|_{L_i^{\prime\prime}})$
is equal to the rotation number of the line field
${\cal N}\cap T_\ast L\subset T_\ast L$
defined on the oriented curve
$\gamma_i$
(in the oriented plane bundle
$T_\ast L$).
This rotation number is equal to
$n_2^\prime -n_1^\prime$,
where
$n_1^\prime$
and
$n_2^\prime$
are the numbers of cusps (points from
$V_2$)
lying on
$\gamma_i$
at which the vector 
$v_2$
points inside
$L_1$
and
$L_2$
respectively.

Adding up the Euler numbers for all pieces of
$L^\prime_i$
and
$L_i^{\prime\prime}$
we get the required formula.
\b
\bigskip

\bigskip
\noindent
{\bf Proof of 
Proposition~\ref{prop-22.3}.}
Obviously, if
$L$
and
${\cal N}$
are oriented, then
$V_1$
subdivides
$L$
into two parts with the common boundary
$V_1$.

Suppose now that
$V_1$
subdivides
$L$
into the two parts
$L_1$
and
$L_2$ 
(``positive'' and ``negative'') as above. 
One can easily define
${\cal N}$ 
near 
$V_1$
so that 
${\cal N}$
has the prescribed singularities of tangency to
$L$
along
$V_1$
and so that for any 
$x$
near
$V_1$,
$x\not\in V_1$,
the oriented planes
$T_x (L)$
and
${\cal N}_x$
together define the positive (resp. the negative) orientation of the 
ambient space
$T_x(T^\ast L)$,
if
$x$
lies in 
$L_1$ (resp.
$L_2$). 
To extend of
${\cal N}$
to the rest of
$L$
means to extend a section of a fiber bundle with a contractible fiber from a 
submanifold of the base to the whole base (see Section~\ref{prelim}
and the proof of 
Proposition~\ref{prop-22.1}). This can always be done.
\b
\bigskip

\subsection{ Proofs of Corollaries~\ref{cor-0.11},\ref{cor-0.1}}
\label{proof-A-II}
\bigskip

First, we use the Lagrangian neighborhood theorem in the same way as before 
but instead of the singularities of the Lagrangian projection
$L\to M$
we consider the singularities of tangency of the zero-section of
$T^\ast L$
to the foliation
${\cal F} = F^\ast (Vert)$,
where
$Vert$
is the foliation of 
$T^\ast M$
by the fibers of the cotangent bundle and
$F$
is a symplectomorphism between a neighborhood of the zero-section of
$T^\ast L$
and a neighborhood of 
$L$
in
$T^\ast M$.
Observe that the Euler class
$e({\cal F})$
is necessarily even.

\bigskip
\noindent
{\bf Proof of Corollary~\ref{cor-0.11}.}
According to Proposition~\ref{prop-22.3},
there exists an oriented 
$\Sigma^2{\hbox{\rm -nonsingular}}$
Lagrangian distribution 
${\cal N}$
along
$L$,
such that the chain on
$L$
associated with 
${\cal N}$
is equivalent to the chain
$\{ T^\ast L\supset L\supset V_1\supset V_2 ; v_1 , v_2\}$. 
Then we finish the proof by applying Theorem A and  
Propositions~\ref{prop-22.1} 
and 
\ref{prop-22.2}
to the distributions
${\cal N}$
and
${\cal F}$.
\b
\bigskip

\bigskip
\noindent
{\bf Proof of Corollary~\ref{cor-0.1}.}
Consider the chain 
${\cal V}=\{V_j, v_j\}$
on the zero-section
$L\subset T^\ast L$
associated with the
$\Sigma^2{\hbox{\rm -nonsingular}}$
oriented Lagrangian distribution
${\cal F}$.
Let
$L_1$
and
$L_2$
as before denote the two parts of
$L$
with the common boundary
$V_1$.

If the submanifold 
$V_2$
is non-empty, we construct another chain on
$L$
in the following way.
We take the chain
${\cal V}$,
choose some points
$x_1,\ldots, x_k$
in the interior part of
$L_1$
and some points
$y_1,\ldots, y_m$
in the interior part of
$L_2$.
Then we add to
$V_1\setminus V_2$
some
$k+m$
new small circles around the points
$x_1,\ldots, x_k, y_1,\ldots, y_m$,
and after that we remove 
$V_2$
from the chain (i.e. we erase all points of
$V_2$
from the curves which form
$V_1$).
In this way we obtain a new chain which we call
${\cal V}^\prime$.
Using the formula from 
Proposition~\ref{prop-22.2}
and according to 
Proposition~\ref{prop-22.3}, 
since
$e({\cal F})$
is even, we get that the number of points in the original submanifold
$V_2$
is even. Therefore  we can find some appropriate numbers
$k$
and
$m$
so that for the new chain
${\cal V}^\prime$
and for the new parts
$L_1^\prime$,
$L_2^\prime$
of
$L$
(defined in the same way as 
$L_1$ 
and
$L_2$
above) we have the following formula:
$\chi (L_1^\prime) - \chi (L_2^\prime) = e({\cal F})$.
We finish the proof by applying Corollary\ref{cor-0.11}
to the chain
${\cal V}^\prime$.
\b

\subsection{ Proof of Theorem A}
\label{proof-A-I}

{\sl
An outline of the proof.
}
As before we  assume that 
$\cal F$
and
$\cal N$
coincide near
$\partial L$.
Our strategy is to reduce Theorem A to Theorem B. 
Initially 
${\cal F}$
and
$\cal N$
may not be homotopic in the class of
$\Sigma^2{\hbox{\rm -nonsingular}}$
oriented Lagrangian distributions along
$L$.
Indeed,
Lagrangian singularities of the type
$\Sigma^2$
are (generically) of codimension 3 and therefore
(generically) at some discrete moments a homotopy between 
${\cal F}$
and
$\cal N$
may pass through a 
$\Sigma^2{\hbox{\rm -nonsingular}}$
distribution.
In order to overcome this problem
we construct some special Hamiltonian isotopy 
$\{ h_t(L)\}$,
$0\leq t\leq 1$,
supported in an arbitrarily small 
neighborhood of a point in
$L$,
so that in the course of this isotopy at exactly one moment
$t$
the Lagrangian submanifold
$h_t (L)$
has a single point of 
$\Sigma^2{\hbox{\rm -singularity}}$ 
of tangency to the fixed foliation
$\cal F$. 
Then we replace
${\cal F}$
by a new Lagrangian distribution
$h^{-1}_1({\cal F})$
along
$L$.
We show that we can apply to
$L$
a number of  isotopies of the same kind as 
$\{ h_t\}$,
each time changing 
$\cal F$
as above, so that at the end instead of the original
$\cal F$
we get a Lagrangian distribution
${\cal F}_1$ 
homotopic to
$\cal N$
in the class of
$\Sigma^2{\hbox{\rm -nonsingular}}$
oriented Lagrangian distributions along
$L$.
Then we apply Theorem B to the distributions
${\cal F}_1$ 
and
${\cal N}$.
This finishes the outline of the proof.
\bigskip

We start with considering the same setup as in the proof of 
Proposition~\ref{prop-22.1}.
Namely, we consider a vector bundle 
$E\to L$
which to each point 
$x\in L$
associates the tangent space
$T_x (T^\ast L)$.
Consider a fiber bundle
$\Lambda_2^{+} (E)\to L$ 
which associates to each point 
$x\in L$
the Grassmannian of oriented Lagrangian planes in the fiber of
$E$
over
$x$.
Both
$\cal F$
and
$\cal N$
can be viewed as sections of
$\Lambda_2^{+} (E)\to L$.
We triangulate
$L$
into small simplices so that over each 2-simplex the vector bundle
$E\to L$
is trivial.

\begin{lem}
\label{lem-2-dim-case-1}
Suppose that the Lagrangian distributions
${\cal F}$
and
${\cal N}$
are homotopic as oriented Lagrangian distributions 
along a closed simplex
$\Delta$
from the triangulation and coincide near
$\partial\Delta$.
Then there exists a Hamiltonian isotopy 
$\{ h_t (\Delta)\}$,
$0\leq t\leq 1$,
such that along
$\Delta$
the Lagrangian distribution
${(h_1^{-1})}_\ast ({\cal F})$
is homotopic to 
${\cal N}$
in the class of
$\Sigma^2{\hbox{\it -nonsingular}}$
Lagrangian distributions (and this homotopy is supposed to be identical near
the boundary of
$\Delta$). 

The Hamiltonian isotopy can be chosen in such a way that it is arbitrarily
$C^0{\hbox{\it -small}}$,
supported inside
$\Delta$
and so that 
$\displaystyle \max_{t,\gamma}
\mid \int_{h_t(\gamma)} pdq \mid$ 
is also arbitrarily
$C^0{\hbox{\it -small}}$.

\end{lem}
\bigskip

\bigskip
\noindent
{\bf Proof of Lemma~\ref{lem-2-dim-case-1}.}
We can assume that the simplex
$\Delta$ 
is sufficiently small. Then, 
without loss of generality, we can identify 
$\Delta$ 
with a neighborhood 
$U\subset {\bf R}^2$
of zero inside the zero-section
${\bf R}^2\subset {\bf R}^4$.
We can also identify the restriction of the vector bundle
$E\to L$
on
$\Delta$
with the tautological vector bundle which associates to each point of
${\bf R}^2\subset {\bf R}^4$
the ambient space
${\bf R}^4$.
Since
$\cal F$
and
$\cal N$
coincide near
$\partial\Delta$
these two distributions along
$\Delta$
together form a spheroid
$S(\Delta)$
in
$\Lambda^{+}_2$.
By the hypothesis of the lemma this spheroid is contractible.
Let us choose a particular homotopy 
$\{ S_t(\Delta)\}$,
$0\leq t\leq 1$,
which contracts the spheroid 
$S(\Delta)$
to a point in
$\Lambda^{+}_2$.
We can assume that this homotopy is generic. Then there is a finite number of 
moments
$t$
at each of which the deformation
$\{ S_t(\Delta)\}$
passes transversally through exactly one of the points
$o^{+}, o^{-}\in\Lambda^{+}_2$.
Without loss of generality we can consider the case when there is exactly one
such moment
$t$
and at this moment the deformation
$\{ S_t(\Delta)\}$
passes transversally through the point
$o^{+}$.

\begin{sublem}
\label{lem-2-dim-case-2}
There exists a Hamiltonian isotopy
$\{ h_t (U)\}$,
$0\leq t\leq 1$,
of
$U$
in
${\bf R}^4$
which satisfies the following property:

\smallskip
\noindent
for any
$t$
the Lagrangian submanifold
$h_t (U)$
coincides with
$U$
near
$\partial U$
so that under the Gaussian map
$h_t (U)$
and
$U$ 
together form a spheroid in
$\Lambda^{+}_2$
which is homotopic to the spheroid
$S_t (\Delta)$
in
$\Lambda^{+}_2\setminus \{ o^{-}, o^{+}\}$.

The Hamiltonian isotopy can be chosen in such a way that it is arbitrarily
$C^0{\hbox{\it -small}}$,
supported inside
$U$
and so that 
$\displaystyle \max_{t,\gamma}
\mid \int_{h_t(\gamma)} pdq \mid$ 
is also arbitrarily
$C^0{\hbox{\it -small}}$.

\end{sublem}

\noindent
{\bf Proof of Sublemma~\ref{lem-2-dim-case-2}.}
Let us pick a point
$x\in \Delta$
at which the planes
$T_x L$
and
${\cal F}_x$
together define a negative 
(with respect to the standard symplectic form)
orientation on the ambient symplectic space
$T^\ast L$
at
$x$.
(If such an
$x$
does not exist, it means that
$\cal F$
and
$\cal N$
are already homotopic as
$\Sigma^2{\hbox{\rm -nonsingular}}$
oriented Lagrangian distributions along
$\Delta$).

We use the Darboux theorem and reduce our problem to a local model in 
${\bf R}^4$,
where
$L$
is identified with the (``horizontal'')
$(q_1,q_2){\hbox{\rm -plane}}$
and
$\cal F$ 
is identified with the foliation of
${\bf R}^4$
by the planes parallel to the ("vertical")
$(p_1,p_2){\hbox{\rm -plane}}$.
Here
$p_1,p_2,q_1,q_2$
denote the standard symplectic coordinates on
${\bf R}^4$.
We identify the (oriented)
$(p_1,p_2){\hbox{\rm -plane}}$
in
${\bf R}^4$
with the point
$o^{+}\in \Lambda_2^{+}$,
so that the (oriented)
$(q_1,q_2){\hbox{\rm -plane}}$
together with
$o^{+}$
define a negative orientation of
${\bf R}^4$
with respect to the standard symplectic form
$dp\wedge dq$.
We will also use the notation
$r= 1/2 (q_1^2 + q_2^2)$,
$R= 1/2 (p_1^2 +p_2^2 +q_1^2 + q_2^2)$.

Given a number
$\delta$
one can find a small enough number
$\epsilon$
and a big enough number
$N$
to construct a family of functions
$\{ f_\tau (q_1, q_2)\}$,
$-\epsilon\leq \tau \geq \epsilon$,
so that the following conditions are fulfilled for any
$\tau$:

\noindent
(i) the function is equal to 
$0$, 
when
$r\geq 3\delta$,
and equal to
$- q_1^2 q_2 + {q_2^3}/3$,
when
$\delta/10\leq r\leq 2\delta$;

\noindent
(ii)
$f_\tau (q_1, q_2) = - q_1^2 q_2 + {q_2^3}/3 + \tau (q_1^2 +q_2^2)$
when
$r\leq \delta/N$;

\noindent
(iii) 
$f_\tau$
does not have critical points when
$\delta/N\leq r\leq \delta/10$.

Consider a Lagrangian submanifold 
$L_\tau\subset {\bf R}^4$
generated by the function
$f_\tau$.
Each
$L_\tau$
has a nonsingular projection on
$L$
and is Hamiltonian isotopic to
$L$
(and this Hamiltonian isotopy is identical at infinity).

Let
$H(R)$
be non-negative function, linear and decreasing on
$[0,\delta/2]$,
decreasing and concave up on
$[\delta/2,\delta]$ 
and equal to zero outside of
$[0,\delta]$.
Consider the Hamiltonian flow 
$\{ h_t^\prime\}$,
$0\leq t\leq T$,
$T= -\pi/2 H_R(0)$,
corresponding to the Hamiltonian function
$H$
on
${\bf R}^4$,
where by
$H_R$
we denote the derivative of
$H$
with respect to
$R$.
Let us apply this Hamiltonian flow to the deformation
$L_\tau$.
A direct calculation shows that for any
$\tau < 0$
and any
$0\leq t\leq T$
the Lagrangian submanifold
$h_t^\prime (L_\tau)$
has no
$\Sigma^2{\hbox{\rm -singularities}}$
of the Lagrangian projection on the zero-section
$L$.
It also shows that for any
$\tau\geq 0$
there exists a unique
$t$,
$0\leq t\leq T$,
such that the Lagrangian submanifold
$h_t^\prime (L_\tau)$
has a single point where its tangent plane coincides with the vertical plane
$o^{+}$
and has no other
$\Sigma^2{\hbox{\rm -singularities}}$
of the Lagrangian projection on
$L$.
One also sees from the calculation that under the Gaussian map
the Lagrangian submanifolds
$h_T^\prime (L_{-\epsilon})$
and
$h_T^\prime (L_{\epsilon})$
together form a small (basic) spheroid around the point
$o^{+}\in \Lambda_2^{+}$
(i.e. this spheroid is homotopic in
$\Lambda^{+}_2\setminus\{ o^{+} , o^{-}\} $
to a boundary of a small ball in
$\Lambda^{+}_2$
centered at
$o^{+}$).
Let us denote by
$\{ h_t\}$
the Hamiltonian isotopy between
$L$
and
$h_T^\prime (L_{\epsilon})$
that we have constructed. 
This Hamiltonian isotopy realizes (under the Gaussian map)
a transversal passage through 
$o^{+}$
in 
$\Lambda_2^{+}$.
To realize a passage in the opposite direction
we need to change simultaneously
$t$
to
$-t$
in the construction of
$h_t^\prime$
and
$\tau$
to
$-\tau$
in the construction of
$\{ f_\tau\}$.

One easily sees that by taking sufficiently 
$C^0{\hbox{\rm -small}}$
deformation
$\{f_\tau\}$
and sufficiently 
$C^1{\hbox{\rm -small}}$
Hamiltonian function
$H(R)$
we can make the Hamiltonian isotopy 
$\{ h_t\}$
arbitrarily 
$C^0{\hbox{\rm -small}}$,
supported inside an arbitrarily small neighborhood of
$x$
and such that 
$\mid\int_{h_t(\gamma)} pdq\}\mid$
is uniformly arbitrarily
$C^0{\hbox{\rm -small}}$
for any
$t$
and any path
$\gamma$
in
$L$. 
This finishes the construction of the special Hamiltonian isotopy
$\{ h_t\}$.
Sublemma~\ref{lem-2-dim-case-2} 
is proved.
\b
\bigskip

\bigskip
\noindent
{\bf End of the proof of Lemma~\ref{lem-2-dim-case-1}.}
By applying a necessary number of times the Hamiltonian isotopy
constructed in Sublemma~\ref{lem-2-dim-case-2}
we can change
$\cal F$,
and accordingly the spheroid
$S(\Delta)$,
so that this spheroid becomes contractible in
$\Lambda^{+}\setminus \{ o^{-}, o^{+}\}$.
It means that we can construct a Hamiltonian isotopy 
$\{ h_t (\Delta)\}$
such that the oriented Lagrangian distribution
${(h_1^{-1})}_\ast ({\cal F})$
is homotopic to 
${\cal N}$
in the class of
$\Sigma^2{\hbox{\rm -nonsingular}}$
Lagrangian distributions along
$\Delta$.
Lemma~\ref{lem-2-dim-case-1} is proved. 
\b
\bigskip

\bigskip
\noindent 
{\bf End of the proof the theorem.}
The obstruction to the homotopy between
$\cal F$
and
$\cal N$
in the class of
$\Sigma^2{\hbox{\rm -nonsingular}}$
Lagrangian distributions along the 1-skeleton of
$L$
lies in
$H^1(L,\partial L; \pi_1(\Lambda^{+}_2\setminus \{ 2\ {\rm points}\}))\cong
H^1(L,\partial L; {\bf Z})$.
As in Proposition~\ref{prop-22.1}
one can see that this obstruction is equal to
zero because
$c({\cal F}) = c({\cal N})$
by the hypothesis of the theorem and Proposition~\ref{prop-22.1}.
Hence we can deform
$\cal F$
into
$\cal N$
in the class of
$\Sigma^2{\hbox{\rm -nonsingular}}$
oriented Lagrangian distributions near the 1-skeleton of
$L$.

By the hypothesis of the theorem and by Proposition~\ref{prop-22.1} we have that
$e({\cal F}) = e({\cal N})$.
Therefore using the standard obstruction theory methods
we can choose the previous homotopy over the 1-skeleton so that
for each closed 2-simplex
$\Delta$
of the triangulation the distributions
${\cal F}$
and
${\cal N}$
coincide near
$\partial\Delta$
and are homotopic as oriented Lagrangian distributions along
$\Delta$.

Then we apply Lemma~\ref{lem-2-dim-case-1}
to each 2-simplex
$\Delta$
and get a new oriented Lagrangian distribution 
${(h_1^{-1})}_\ast ({\cal F})$
along
$L$
which coincides with
$\cal N$
near
$\partial L$
and which is homotopic to
$\cal N$
in the class of
$\Sigma^2{\hbox{\it -nonsingular}}$
Lagrangian distributions along the whole 
$L$.
Therefore Theorem B can be applied and this finishes the proof of Theorem A.
\b
\bigskip

\section{Singularities and generating functions}
\label{gen-f}

Let us now consider the case of
an immersed Lagrangian submanifold
$L\subset T^\ast M$,
which is generated by a function 
$f$
defined on an open subset
$U$ 
of
$M\times K^N$,
where
$K^N$
is a (connected) manifold of dimension
$N$.

\begin{defin}
\label{def-1.5}
{\rm 
In a situation as above let us define
$V_L\subset U\subset M\times\KN$ 
as   
$V_L=\{ (q,\ksi)\ |\ {\partial f}/{\partial\ksi}(q,\ksi) =0\}$ 
(by 
$(q,\ksi)$
we denote local coordinates on
$M\times K^N$).
}
\end{defin}
\bigskip

There exists a natural diffeomorphism
$F$ 
between
$V_L$ 
and
$L$.
The diffeomorphism
$F$ 
maps  the singularities of the natural projections  
$P:V_L\to M$
and 
$\pi:L\to M$
onto each other. 
On the other hand, as it follows from Section~\ref{ch-dirsurg}, 
we can use Definition~\ref{def-1.4} 
and construct a chain
${\cal V}$
on
$L$.
Under the diffeomorphism
$F$
the chain
${\cal V}$
is mapped into a chain
${\cal V}_1$
on
$V_L$. 

\begin{defin}
\label{def-1.5a}
{\rm
The chain 
${\cal V}_1$
on
$V_L$
defined as above is called the
{\it
chain induced by
$L$.
}
}
\end{defin}
\bigskip

Let us now consider Case II mentioned in Section~\ref{ch-dirsurg}, 
i.e. the case when
${\rm dim}\ K = 1$
and
$L$
is generated by a function defined on (an open subset of)
$M\times K^1$.
In this case in addition to the chain
${\cal V}_1$
induced on
$V$
by
$L$
we can define a similar object called
$\alpha{\hbox{\it -chain}}$ 
which, unlike 
${\cal V}_1$,
uses some information about
the function
$f$.

\begin{defin}
\label{def-1.6}
{\rm 
\cite{El}
Let
$V^n\subset M^n\times K^1$
be a submanifold.
An
$\alpha{\hbox{\it -chain}}$
on
$V$
is a weak chain 
${\cal V}^\prime =\{V_i, v_j^\prime \}, i\geq {-1},j\geq 2,$
on
$V$
equipped with a unit vector field
$v^\prime_1$
defined over
$V_1\setminus V_2$.
This vector field 
$v^\prime_1$
has to satisfy the following properties:

\smallskip
\noindent
(i) 
$v^\prime_1$
is normal to
$V_0=V$
in
$V_{-1}=M\times K^1$;

\smallskip
\noindent
(ii)
$v^\prime_1$
cannot be extended
(as a vector field normal to 
$V_0$
in
$V_{-1}$) 
to any subset 
$C\subset V_1$ 
that has a nontrivial  intersection with 
$V_2$.
}
\end{defin}

\begin{rem}
\label{rem-1.6}
{\rm
In \cite{El}
$\alpha{\hbox{\rm -chains}}$
are called simply 
{\it
chains.
}
}
\end{rem}

\begin{defin}
\label{def-1.3a}
{\rm
Two 
$\alpha{\hbox{\rm -chains}}$ 
on two submanifolds
$V_1, V_2\subset M\times K^1$
are called
{\it equivalent,}
if there exists an isotopy between
$V_1$
and
$V_2$
in
$M\times K^1$
which maps one
$\alpha{\hbox{\rm -chain}}$
into another.
}
\end{defin}

\bigskip
\noindent
{\bf Notational convention.}

\smallskip
\noindent
By an (oriented)
{\it
line field along a submanifold
$V\subset M\times K^1$
}
we mean a field of (oriented) lines 
in
$T_\ast (M\times K^1)$
defined along
$V$.
\bigskip

\begin{defin}
\label{def-1.6a}
{\rm \cite{El}
Let
$V\subset M\times K^1$
be a submanifold, and let
${\cal L}$
be an (oriented)  line field along
$V$.
We define an
{\it
$\alpha{\hbox{\it -chain}}$
on
$V$
associated with
${\cal L}$
}
in the following way. First, we define a 
{\it
weak chain on
$V$
associated with
${\cal L}$
}
in the same way as we did in Definition~\ref{def-1.4}.
Then we define the vector field
$v_1^\prime$
in the same way as we defined the vector fields
$v_j$,
$j\geq 2$,
in Definition~\ref{def-1.4}. 
Namely, at each point
$x\in V_1\setminus V_2$
the vector
$v_1^\prime (x)$
is a unit normal vector to
$V$
in
$M\times K^1$
such that its projection 
$P_\ast (v_1^\prime (x))$
on
$M$
is an exterior normal to the boundary of
$P(V)\subset M$
at
$P(x)$,
where
$P:M\times K^1\to M$
is the canonical projection.

One can check that since
${\rm dim}\ K=1$
this construction indeed defines an
$\alpha{\hbox{\rm -chain}}$ 
and this
$\alpha{\hbox{\rm -chain}}$
is unique up to the equivalence of
$\alpha{\hbox{\rm -chains}}$.
}
\end{defin}
\bigskip

Let
$Vert$
be a line field along
$V_L$
which to each point
$x\in V_L$
associates the tangent space to the fiber of the projection
$M\times K^1\to M$
at
$x$.

\begin{defin}
\label{def-1.6b}
{\rm
If
$L$
is generated by a function
$f$,
then the
{\it
$\alpha{\hbox{\it -chain}}$
on
$V_L$
generated by 
$f$
}
is an
$\alpha{\hbox{\it -chain}}$
on
$V_L$
associated with the line field 
$Vert$.
}
\end{defin}

\begin{rem}
\label{rem-1.2a}
{\rm
Observe that the line field
$Vert$
is tangent to
$V_i$ 
along
$V_{i+1}$
and transversal to
$V_i$
along
$V_i\setminus V_{i+1}$,
$i=1,\ldots ,n$.
At any point
$x\in V_1\setminus V_2$
the vector
$v_1 (x)$ 
from the chain on
$V$
induced by
$L$
determines an orientation of
$Vert$
at
$x$.
It means that the first two derivatives of
$f$
along 
$Vert$
(in the direction given by the orientation)
vanish at
$x$
and the third one is always positive.
Similarly, if
$x\in V_i\setminus V_{i+1}$,
$i>1$,
then the first 
$i+1$ 
derivatives of
$f$
along the oriented line field
$Vert$
vanish at 
$x$
and the
$(i+2){\hbox{\rm -nd}}$
one is nonzero.

We can also describe the 
$\alpha{\hbox{\rm -chain}}$ 
on
$V_L$
generated by 
$f$
in the following way.
A choice of orientation of
$K^1$
makes
$Vert$ 
globally oriented.
This allows us to  define correctly the sign of the
$(i+2){\hbox{\rm -nd}}$
derivative of 
$f$
along the positive direction of 
$Vert$
at each point
$x\in V_i\setminus V_{i+1}$,
$i\geq 1$.
We call this sign
$sign_1 (x)$.

On the other hand, if we fix an orientation of 
$K^1$
then the signs of the partial derivatives of
$f$
with respect to
$\xi$
simultaneously coorient
$V_i$ 
in
$V_{i-1}$,
$i\geq 0$.
Therefore to each point 
$x\in V_i\setminus V_{i+1}$,
$i=1,\ldots, n$,
we can associate another  sign called
$sign_2 (x)$ 
which is taken to be plus, if the vector 
$v_i^\prime (x)$ 
from the 
$\alpha{\hbox{\rm -chain}}$ 
generated by
$f$
agrees with the previous coorientation of
$V_{i-1}$
in
$V_{i-2}$,
and minus otherwise.
By a suitable choice of orientation of
$K^1$
we can get that
$sign_1(m) = sign_2 (m)$
for any 
$x\in V_1$.

Thus instead of vector fields
$v_i^\prime$
in the 
$\alpha{\hbox{\rm -chain}}$ 
generated by 
$f$
we can work with the signs of the derivatives 
${{\partial}^i f}/{\partial {\xi}^i}$.
}
\end{rem}

\begin{figure}
\centerline{\psfig{figure=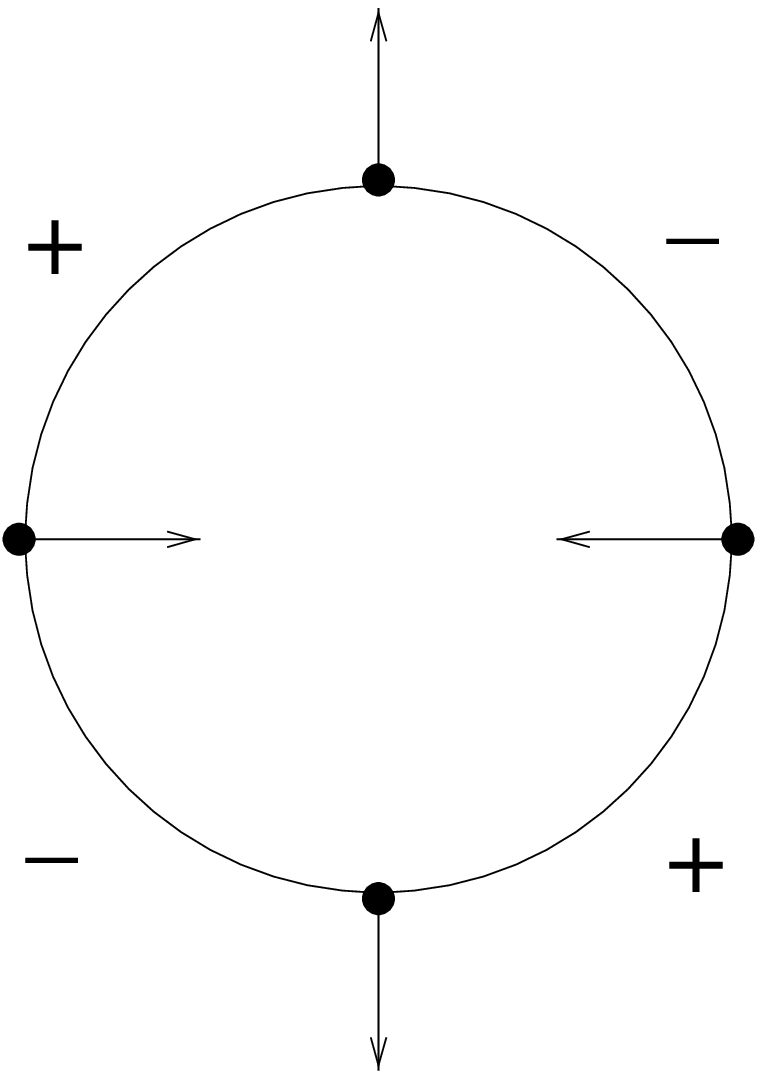,height=80mm}}
\caption{The 
$\alpha{\hbox{\rm -chain}}$ 
corresponding to the Lagrangian collapse.}
\label{fig5}
\end{figure}

\begin{exam}[``Lagrangian collapse'']
\label{exam-1.1a}
{\rm 
Consider a Lagrangian submanifold
$L_\varepsilon\subset  T^\ast ({\R}^2)$
generated by a function
\[f_\varepsilon (q_1 , q_2 ,\xi)=q_1 cos \xi +q_2 sin \xi - 
\chi_\varepsilon (r) cos 2\xi ,\]
where
$r = q_1^2 +q_2^2$
and
$\chi_\varepsilon (r)$
is a function such that:

\smallskip
\noindent
(i) 
$0\leq \chi_\varepsilon (r)\leq \varepsilon/2$
for any
$r$;

\smallskip
\noindent
(ii)
$\chi_\varepsilon (r) = \varepsilon/2$,
if
$0\leq r\leq 100$;

\smallskip
\noindent
(iii)
$\chi_\varepsilon (r) = 0$,
if
$r\geq 1000$.

\smallskip
For a sufficiently small
$\varepsilon$
we get that
$L_\varepsilon$
is an embedded Lagrangian submanifold and
${\partial}^2 f/{\partial {\xi}^2} (q,\xi) \neq 0$, 
if
$q_1^2 + q_2^2 > 100$. 

The chain on the corresponding cylinder
$V_{L_\varepsilon}$
induced by
$L_\varepsilon$
is shown on 
Fig.~\ref{fig1}.
It consists of a circle
of folds 
($V_1$) 
with four cusps on it 
($V_2$). 
The directions of the vector field
$v_1$
on
$V_1$
are shown by the small arrows, and the directions of
$v_2$
on
$V_2$
are shown by the big arrows.

The 
$\alpha{\hbox{\rm -chain}}$ 
on
$V_{L_\varepsilon}$
generated by
$f_\varepsilon$
is shown on 
Fig.~\ref{fig5}.
The signs attached to each connected component of
$V_1\setminus V_2$
on the pictures are the signs of 
${{\partial}^3 f_\varepsilon }/{\partial {\xi}^3}$
on those components
(see Remark~\ref{rem-1.2a}). 
}
\end{exam}
\bigskip

To define a direct surgery for
$\alpha{\hbox{\rm -chains}}$
we need to adjust the Definitions~\ref{def-1.7} 
and \ref{def-1.8}.

\begin{defin}[Adjusted Definition~\ref{def-1.7} 
for
$\alpha{\hbox{\bf -chains}}$]
\label{def-1.7a}
{\rm
A
{\it
basis 
of a} 
{\it
direct
surgery
}
of order
$s\geq 1$
and index
$p\geq 0$
is formed by an
$\alpha{\hbox{\rm -chain}}$
${\cal V}^\prime$
on
$V_L$
and an embedding
$\phi :D^p\times [-1,1]\to V_{s-1}$.
This data must satisfy the same conditions as in
Definition~\ref{def-1.7}
except for the following: 

\smallskip
\noindent
(i) there is no need for the vectors
$\nu_3$,
$\nu_4$,
when
$p=0$,
$s=1$;

\noindent
(ii) the condition e) is not needed;

\noindent
(iii) in the condition d) one should write
$s\geq 1$
instead of
$s>1$.
}
\end{defin}

\begin{defin}[Adjusted Definition~\ref{def-1.8} for
$\alpha{\hbox{\bf -chains}}$]
\label{def-1.8a}
{\rm
Suppose
}
\break
{\rm
that an 
$\alpha{\hbox{\rm -chain}}$
$\cv$ 
on
$V_L$
and an embedding
$\phi :D^p\times [-1,1]\to V_{s-1}$
(possibly with a pair of vectors
$\nu_1 ,\nu_2$)
form a basis of a direct surgery of order
$s$
and index 
$p$.
The result of the direct surgery associated with the basis
is a new chain
${\cal V}^\prime$
on
$V_L$
defined in the same way as in Definition~\ref{def-1.8}
but we should not single out the case
$p=0$,
$s=1$ 
(i.e. the same construction applies to any
$p$
and
$s$).
}
\end{defin}

\begin{rem}
\label{rem-1.3}
{\rm
Suppose that
$L$
is generated by a function 
$f$.
A basis of a direct surgery for the chain 
on
$L$
induces a corresponding basis of a direct surgery
for the chain on
$V_L$
induced by
$L$.

One can check that any basis of a direct surgery for the 
$\alpha{\hbox{\rm -chain}}$ 
on
$V_L$,
generated by
$f$,
induces a basis of a direct surgery for the chain on
$V_L$
induced by
$L$,
{\sl
but not necessarily vice versa!
}

Indeed, if we want to realize a direct surgery  of order 
$s$
along a disk
$D^p$
for the 
$\alpha{\hbox{\rm -chain}}$ 
generated by
$f$,
we have to have that the signs of 
${{\partial}^{s+2} f}/{\partial {\xi}^{s+2}}$
extend from the boundaries of each of the disks
$D^p\times \{ \pm 1\}$
to the disk itself. This might not hold for a basis of a direct 
surgery for the
chain on
$V_L$
induced by
$L$,
if
$p=1$.
}
\end{rem}
\bigskip

As before we can introduce the notion of an 
{\it
inverse surgery
}
for
$\alpha{\hbox{\rm -chains}}$
as an operation inverse to a direct surgery.
For an accurate definition of an inverse surgery as a
formal operation on
$\alpha{\hbox{\rm -chains}}$
see \cite{El}.
The analogs of Propositions~\ref{prop-1.1} 
and 
\ref{prop-1.2} 
for
$\alpha{\hbox{\rm -chains}}$
look as follows.

\begin{prop}
\label{prop-1.1a}
{\rm
\cite{El}
}
Suppose that we have a generic 1-parameter family  
$\{ {\cal L}^t\}$
of line fields along
$V$.
Then the 
$\alpha{\hbox{\it -chain}}$ 
on
$V$
associated with
${\cal N}^1$
can be obtained from the 
$\alpha{\hbox{\it -chain}}$
associated with
${\cal N}^0$
by a sequence of direct and inverse surgeries.
\end{prop}

\begin{prop}
\label{prop-1.2a}
{\rm
\cite{El}
}
Consider an inverse surgery on an 
$\alpha{\hbox{\it -chain}}$ 
${\cal V} = \{ V_i,v_i\}$
which transforms 
${\cal V}$
into another 
$\alpha{\hbox{\it -chain}}$
${\cal V}^\prime = \{ V_i^\prime, v_i^\prime\}$.
Suppose that the inverse surgery is an inverse of  a direct surgery 
of order
$s$
along a disk
$D^p$
(i.e. the direct surgery transforms
${\cal V}^\prime$
into
$\cal V$).
Suppose also that
there exists an embedding 
$\alpha : D^1\to V_0^\prime$
which satisfies the following properties:

\noindent
(i) 
$\alpha (D^1)\cap V_1^\prime = \alpha (\partial D^1)$
and 
$\alpha (D^1)$ 
is normal to
$V_1^\prime$
at 
$\alpha (\partial D^1)$;

\noindent
(ii)
the vector field 
$v_1$
cannot be extended along
$\alpha (D^1)$
as a normal vector field to
$V_0^\prime$
in 
$V_{-1}^\prime$.
\smallskip

Then the inverse surgery can be represented as a
composition of direct surgeries.

\end{prop}

\begin{defin}
\label{def-1.9a}
{\rm
\cite{El}.
We call an inverse surgery satisfying the hypothesis of Proposition~\ref{prop-1.2a}
{\it 
complete.
}
}
\end{defin}

\section{Proofs of Propositions~\ref{prop-2.1},\ref{thmC}. 
Lagrangian Collapse}
\label{sect-gen-f-pfs}

\bigskip
\noindent
{\bf Notational convention.}

\smallskip
\noindent
From now on by
$B_r^l(x)$
we denote an
$l{\hbox{\rm -dimensional}}$
ball
of radius
$r$
with the center at
$x$.
By
$B_r^l$
we denote an
$l{\hbox{\rm -dimensional}}$
ball of radius
$r$
with the center at 
$0\in {\bf R}^l$.
We omit the lower index for unit balls. By
$S(x,r)$
we denote a sphere of radius
$r$
with the center at
$x$. 
\bigskip

\noindent
{\bf Proof of 
Proposition~\ref{prop-2.1}.}
Let us first sketch the general idea of the proof.
In a neighborhood of
$D^p$
in
$T^\ast L$
we introduce new symplectic coordinates so that
the singularities of tangency of
$L$
to
${\cal F}$
near the disk
$D^p$
become the singularities of Lagrangian projection of a Lagrangian
submanifold in a cotangent bundle. This Lagrangian submanifold 
(in a neighborhood of
$D^p$)
is defined by a generating function.
Then we reduce 
Proposition~\ref{prop-2.1} to a problem concerning a
deformation of this generating function. The latter problem is dealt 
with in the following proposition.

\begin{prop}
\label{prop-3.1}
Let
$L\subset T^\ast M$ 
be an embedded Lagrangian submanifold generated by
a function 
$f$
defined on an open subset of
$\MKN$.
Let
$\cal V$ 
be the chain  
induced by
$L$
on 
$V_L$.
Let 
${\cal V}_1$ 
be the result of a direct surgery on
$\cal V$
of order 
$s$
($s\geq 1$)
along a disk 
$D^p\subset V_{s-1}$
whose projection from
$\MKN$
to
$M$
is an embedded disk in
$M$. 
Then for any sufficiently small neighborhood 
$U$ 
of 
$D^p$ 
in
$\MKN$
there exists a 1-parametric family of functions 
$\{ f_\tau\}$,
$f=f_0$,
$0\leq\tau\leq 1$, 
such that

\smallskip
\noindent
(i) for any
$\tau$,
$0\leq\tau\leq 1$,
the function
$f_\tau$
generates an embedded Lagrangian submanifold
$L_\tau \subset {T^\ast}M$
which is Hamiltonian isotopic to
$L$;

\noindent
(ii) the chain 
${\cal V}_1$ 
is equivalent to the chain
on
$V_{L_1}\cong V_L$
induced by
$L_1$.
\smallskip

The deformation 
$\{ f_\tau\}$ 
can be made arbitrarily
$C^1{\hbox{\it -small}}$
and supported inside
$U$.

If
${\rm dim}\ K = 1$,
${\cal V}^\prime$
is the 
$\alpha{\hbox{\it -chain}}$
generated by
$f$
and
$\alpha{\hbox{\it -chain}}$
${\cal V}_1^\prime$
is the result of a direct surgery on
${\cal V}^\prime$
of order 
$s$
($s\geq 1$) 
along a disk 
$D^p\subset V_{s-1}$,
then we can achieve that
${\cal V}_1^\prime$
is equivalent to the 
$\alpha{\hbox{\it -chain}}$
generated by the function
$f_1$
on
$V_{L_1}\cong V_L$.
\end{prop}
\bigskip

We will prove Proposition~\ref{prop-3.1}
in Section~\ref{pf-prop-3.1-chap}.

\bigskip
\noindent
{\bf End of the proof of 
Proposition~\ref{prop-2.1}.}
Since
$p\leq n-1$,
we can assume that
$D^p$
is transversal to the foliation
${\cal F}$.
Using the Lagrangian neighborhood theorem
(see \cite{We1},\cite{W}) 
we can find a symplectomorphism which identifies
a small open tubular neighborhood 
$U$
of 
$D^p$
in
$T^\ast L$
with a product
${U^\prime}\times B^n_{\epsilon}\subset T^\ast {U^\prime}\cong {\bf R}^{2n}$
for some small
$\epsilon$,
where
${U^\prime}\supset D^p$
is a small tubular neighborhood of
$D^p$
in the Lagrangian subspace
${\bf R}^n\times 0 \subset {\bf R}^{2n}$
and the foliation
${\cal F}$
restricted on
$U$
is formed by the leaves of the Lagrangian projection 
${\bf R}^{2n}\to {\bf R}^n$
restricted on
$U$.
Thus we can present the singularities of tangency of
$L$
to
${\cal F}$
near
$D^p$
as singularities of a Lagrangian projection.

Let
$L^\prime$ 
be the part of
$L$
lying inside
$U$
(and
$U$
itself now is supposed to lie inside
${\bf R}^{2n}$).
Since
$L^\prime$
is contractible 
it follows from a theorem 
of Giroux \cite{G} that
$L^\prime$ 
can be generated by a function
$f$
defined on an open subset of
${U^\prime}\times {\bf R}^N$ 
for some 
$N$.
Since
$D\subset {U^\prime}$
we can apply 
Proposition~\ref{prop-3.1} to the Lagrangian submanifold
$L^\prime\subset T^\ast {U^\prime}$
and the generating function
$f$.
This finishes the proof of 
Proposition~\ref{prop-2.1}.
\b
\bigskip

\begin{figure}
\centerline{\psfig{figure=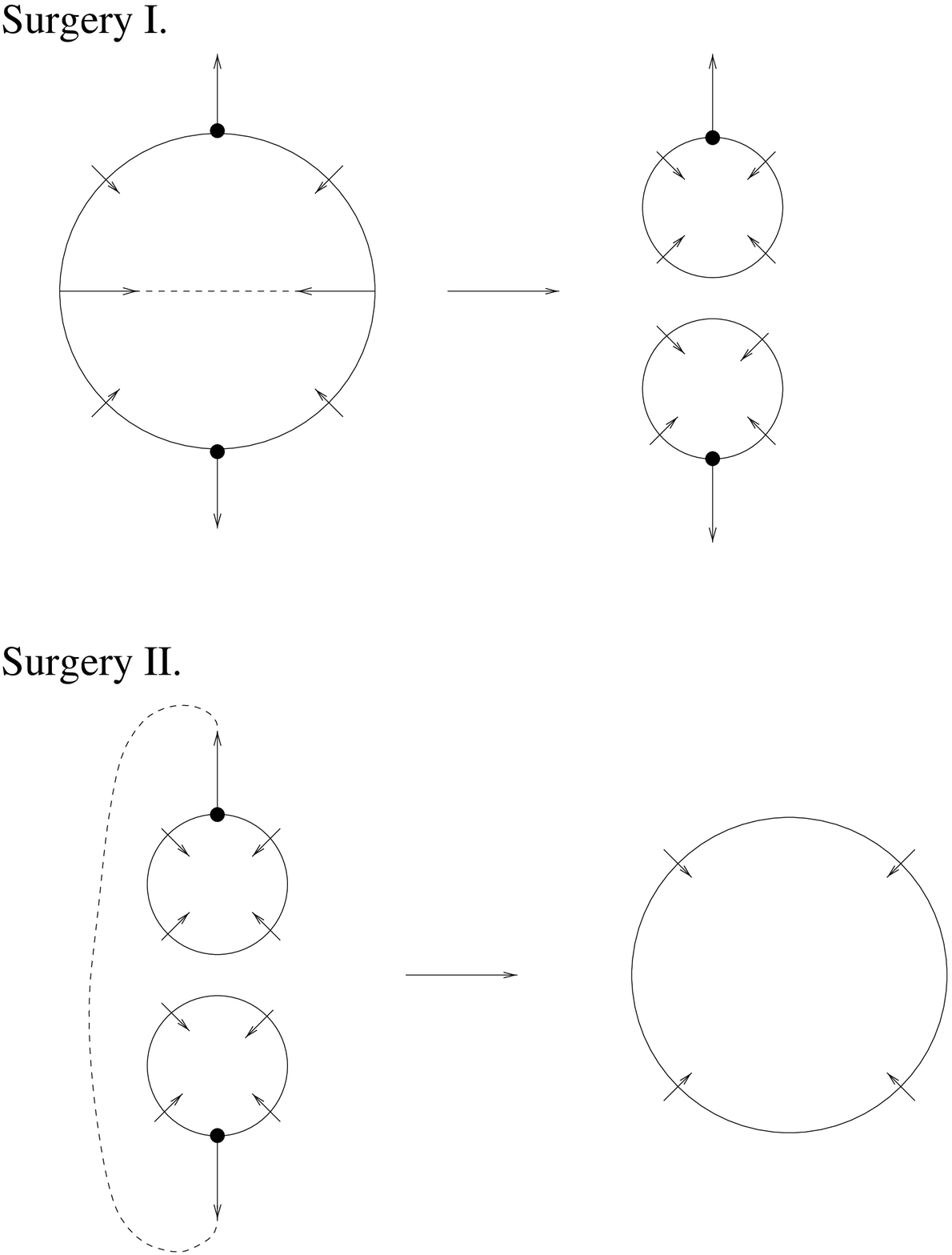,height=175mm}}
\caption{The direct surgeries killing the cusps (but not realizable by a 
deformation of the 
{\it
global 
}
generating function) for the Lagrangian collapse.}
\label{fig6}
\end{figure}

\begin{figure}
\centerline{\psfig{figure=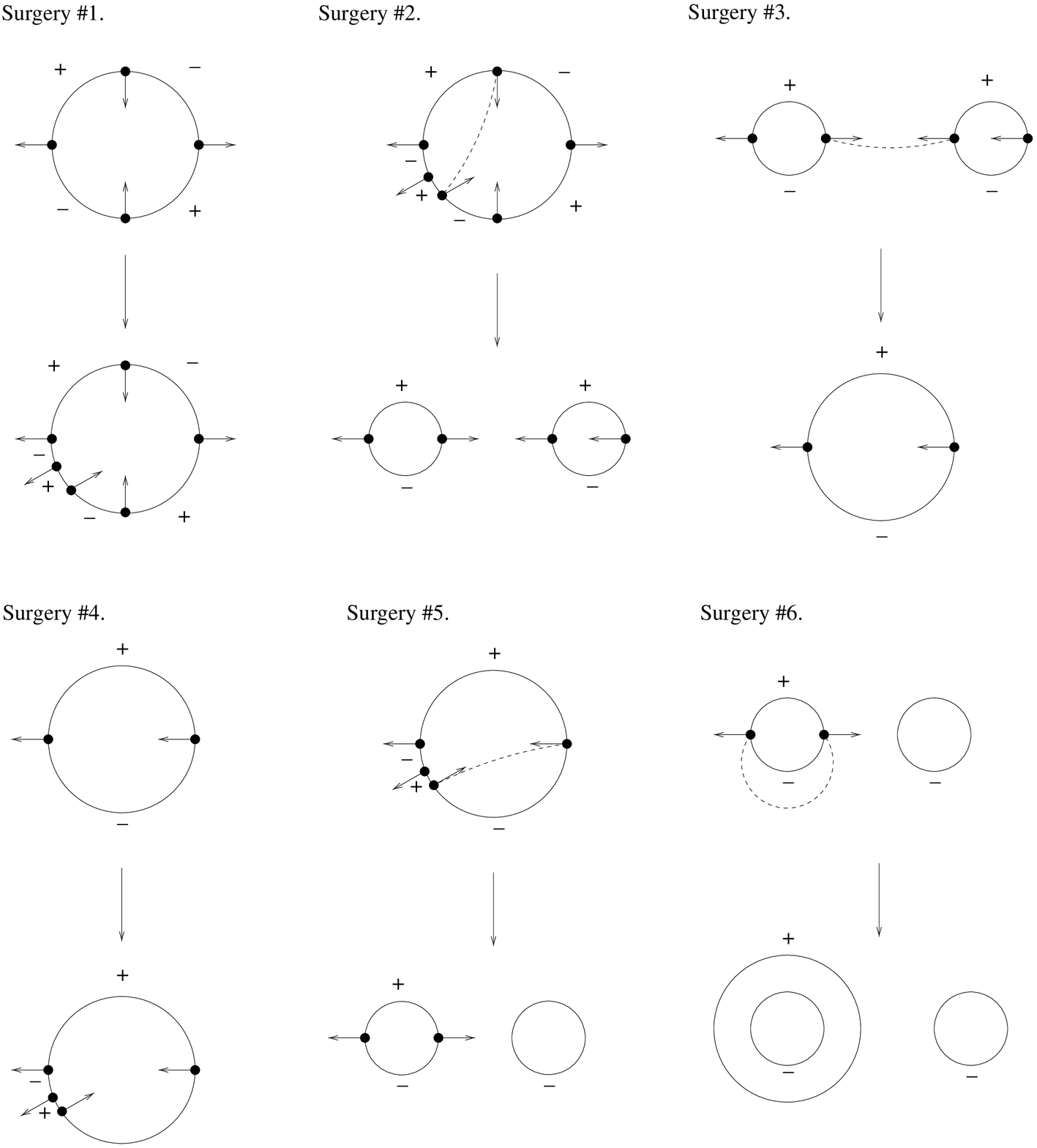,height=175mm}}
\caption{The direct surgeries killing the cusps and realizable by a 
deformation of the
{\it
global 
}
generating function in the case of Lagrangian collapse.}
\label{fig7}
\end{figure}

\begin{exam}
\label{exam-2.1}
{\rm
Consider an embedded Lagrangian submanifold 
$L$
in
$T^\ast {\R}^2$
generated by the function
$f(q_1 , q_2 ,\xi)=q_1 cos \xi +q_2 sin \xi$
defined on
${\R}^2\times S^1$
(as in Examples~\ref{exam-0.1},\ref{exam-1.1a}).
Consider first a deformation of 
$L$
given by the following deformation of
$f$:
\[f_\varepsilon (q_1 , q_2 ,\xi)=q_1 cos \xi +q_2 sin \xi - 
\chi_\varepsilon (r) cos 2\xi\]
(as in Example~\ref{exam-1.1a}).
For 
$\varepsilon$
close enough to 0 the function
$f_\varepsilon$
generates an embedded Lagrangian submanifold 
$L_\varepsilon$
of
$T^\ast {\R}^2$ 
which is Hamiltonian isotopic to 
$L$
(as in Example~\ref{exam-1.1a}).
By choosing a sufficiently small
$\varepsilon$ 
the Hamiltonian isotopy can be made 
compactly supported and arbitrarily
$C^0{\hbox{\rm -small}}$.

We take such an 
$\varepsilon$,
consider the chain on
$V_{L_\varepsilon}$
induced by 
$L_\varepsilon$
(see Example~\ref{exam-1.1} and 
Fig.~\ref{fig1}) 
and apply to it two subsequent direct surgeries as shown on 
Fig.~\ref{fig6}.
Punctured lines indicate the 1-dimensional disks
along which we perform direct surgeries of index 1.

After applying the two surgeries we get  a chain which is formed by
one circle of folds with no cusps on it,
i.e. after the surgery the submanifold
$V_1$
is a circle and
$V_2$
is empty. According to 
Proposition~\ref{prop-2.1}, 
such a chain corresponds to an 
{\it embedded} 
Lagrangian submanifold 
$L_1$
of
$T^\ast {\R}^2$ 
which is Hamiltonian isotopic to 
$L$
and which projects on
${\bf R}^2$
with only fold-type singularities.

However, the two surgeries we made (and hence the whole isotopy) 
cannot be realized by a deformation of the generating function
$f$.
To construct a Hamiltonian isotopy covered by a deformation of
$f$
we have to make direct surgeries for the 
$\alpha{\hbox{\rm -chain}}$ 
generated by 
$f_\varepsilon$ 
(see Example~\ref{exam-1.1a} 
and 
Fig.~\ref{fig5}).
These surgeries are shown on the 
Fig.~\ref{fig7}.
According to 
Proposition~\ref{prop-3.1}, 
we can realize these surgeries
by arbitrarily
$C^0{\hbox{\rm -small}}$
Hamiltonian isotopies of
$L$
covered by arbitrarily 
$C^1{\hbox{\rm -small}}$
deformations of the generating function.
}
\end{exam}

\bigskip
\noindent
{\bf Proof of 
Proposition~\ref{thmC}.} 
First, we need an analogue of
Proposition~\ref{prop-2.2} for
$\alpha{\hbox{\rm -chains}}$.
Similarly to what we did before, we call a pair of 
$(n-1){\hbox{\rm -dimensional}}$
embedded spheres 
$S_1, S_2\subset V_L$
a
{\it
double fold pair of spheres},
if
$S_2^{n-1}$
bounds a ball in
$L$
and the sphere
$S_1^{n-1}$
lies inside that ball.

Consider the 
$\alpha{\hbox{\it -chain}}$
${\cal V}$
on 
$V_L\subset M\times K^1$
generated by
$f$.

\begin{prop}
\label{prop-2.2a}
Let
$B\subset L$
be a ball in
$V_L$
such that the restriction of the projection
$P:V_L\to M$
on
$B$
is nonsingular.
Consider a new 
$\alpha{\hbox{\it -chain}}$
${\cal V}_1$
which is a union of
${\cal V}$
and a double fold pair of spheres lying in
$B$.
Then there exists a deformation 
$\{ f_t\}$,
$0\leq t\leq 1$,
of
$f=f_0$
in the class of generating functions defined on (an open subset of)
$M\times K^1$,
such that 

\smallskip
\noindent
(i) the deformation 
$\{ f_t\}$
covers  a Hamiltonian isotopy of
$L$;

\noindent
(ii) the
$\alpha{\hbox{\it -chain}}$
${\cal V}_1$
is equivalent to the 
$\alpha{\hbox{\it -chain}}$
generated by
$f_1$.
\smallskip

The deformation 
$\{ f_t\}$
can be made arbitrarily
$C^1{\hbox{\it -small}}$
and supported inside 
$B$.

Any of the two possible combinations of the orientations of
$v_1^\prime$ 
on
$S_1$
and
$S_2$
can be realized by some deformation 
$\{ f_t\}$
as above.

\end{prop}
\bigskip

We will prove 
Proposition~\ref{prop-2.2a} 
along with 
Proposition~\ref{prop-2.2} 
in 
Section~\ref{double-folds-pfs}.

\bigskip
\noindent
{\bf End of the proof of 
Proposition~\ref{thmC}}.
Consider the line field
$Vert$
along
$V_L$
(see Definition~\ref{def-1.6b})
and let 
$Norm$
be the normal (oriented) line field along
$V_L$.

The only obstruction for the line fields
$Vert$
and
$Norm$
to be homotopic
as
line fields along
$V_L$
lies in the homotopy group
$\pi_n ({RP}^n)\cong {\bf Z}$.
If 
$L$
is noncompact or has a boundary, we may assume as before that
$Vert$
and
$Norm$
coincide near the boundary and outside of some compact subset of
$V_L$
and we look for a homotopy identical at those parts of
$V_L$.

Suppose that the obstruction is equal to
$k\neq 0$.
Then 
we get rid of the obstruction by creating additional 
double folds on
$L$
and by cutting out  small neighborhoods of the balls bounded by the 
interior spheres in each double fold pair of spheres. 
On the resulting manifold the obstruction vanishes.

Formally speaking, we do the following. We take 
$k$
points
$x_1,\ldots,x_k\in V_L$
which are nonsingular for the  projection
$P:V_L\to M$.
Using 
Proposition~\ref{prop-2.2a} 
we find a
$C^1{\hbox{\rm -small}}$ 
deformation
$\{ f_\tau\}$
of
$f$
which creates additional double folds along
each pair of arbitrarily small concentric spheres 
$S(x_i,\epsilon)$,
$S(x_i,2\epsilon)$,
$i=1,\ldots,k$.
We still denote our manifold by 
$V_L$
and our generating function by
$f$. 

According to 
Proposition~\ref{prop-2.2a}, 
we can arrange the orientations of
$v_1^\prime$ 
on 
$S(x_i,\epsilon)$,
$S(x_i,2\epsilon)$,
$i=1,\ldots,k$,
in such a way that
if we delete all the balls
$B^n_{3/2\epsilon}(x_i)$,
$i=1,\ldots,k$,
from 
$V_L$,
then on the new manifold
$V_L^\prime$
the obstruction vanishes.
Hence
$Vert$
is homotopic to
$Norm$
as a line field 
along the new manifold
$V_L^\prime$.
Let us denote by
$\{ {\cal L}^t\}$
the homotopy between
$Vert$
and
$Norm$
that we have obtained.
We apply 
Proposition~\ref{prop-1.1a} 
to the homotopy
$\{ {\cal L}^t\}$ 
and get a sequence of direct and inverse surgeries,
corresponding to
$\{ {\cal L}^t\}$. 
After that we realize these surgeries by Hamiltonian isotopies as follows.

If we need to realize a 
{\sl
direct surgery
}
we proceed according to 
Proposition~\ref{prop-3.1}.

If at a given step we need to realize
{\sl
a complete inverse surgery}
we use 
Proposition~\ref{prop-1.2a}
and decompose it into a product of direct surgeries. 
Thus the problem is reduced to the previous case.

If at a given step we have an inverse surgery 
which is 
{\sl
not complete, 
}
we make it complete using 
the same trick we used in the proof of Theorem B. 
Namely, we use
Proposition~\ref{prop-2.2a}
and create an appropriate double fold in a neighborhood of an
appropriate point by a
$C^1{\hbox{\rm -small}}$ 
deformation of the generating function.
Then we form a new 
$\alpha{\hbox{\rm -chain}}$ 
by adding the exterior sphere from the double fold pair of spheres
to the chain which is the result of the inverse surgery. We can assume that 
this new chain is the result of the inverse surgery. 
It makes the inverse surgery complete and we reduce the problem to
the previous case.

Thus we have  constructed a sequence of
arbitrarily 
$C^0{\hbox{\rm -small}}$
Hamiltonian isotopies. The composition of all these 
isotopies gives us the required
Hamiltonian isotopy of the original Lagrangian submanifold
$L$.
At each step the deformations of 
the Lagrangian submanifold were covered by a deformation
of generating functions defined on (an open subset of)
$M\times K^1$.
Therefore the final Hamiltonian isotopy of
$L$
also satisfies this property.
\b
\bigskip

By exactly the same argument we can prove the following strengthened version
of 
Proposition~\ref{thmC}.

Suppose again that an embedded  Lagrangian submanifold
$L\subset T^\ast M$
is generated by a function
$f$
defined on an open subset of 
$M\times K^1$,
where
$K$
is a (connected) 1-dimensional manifold.
Let
${\cal N}$
be a generic oriented line field along
$V_L$
homotopic to
$Vert$
as an oriented line field along
$V_L$.
Let
${\cal V}$
be the 
$\alpha{\hbox{\it -chain}}$
on
$V_L$
associated with
${\cal N}$.

\begin{prop}
\label{thmC-full-str}
One can deform 
$L$
by a Hamiltonian isotopy into a new Lagrangian submanifold 
$L_1$
such that the weak chain on 
$V_{L_1}$
induced by
$L_1$
is equivalent to the union of the
$\alpha{\hbox{\it -chain}}$
${\cal V}$
and a number of additional double fold pairs of spheres in
$L$.
This Hamiltonian isotopy can be covered by a  
deformation 
$\{f_t\}$
of the generating function
$f$ 
in the class of generating functions
defined on (an open subset of)
$M\times K^1$.

The deformation 
$\{ f_t\}$
can be made compactly supported and arbitrarily
$C^1{\hbox{\it -small}}$.

\end{prop}
\bigskip

The analogs of
Proposition~\ref{thmC-full-str}
can be proved in the case when
$L\subset T^\ast M$
is an immersed Lagrangian submanifold and in the Legendrian case.

\section{Proof of 
Proposition~\ref{prop-3.1}}
\label{pf-prop-3.1-chap}

We present a proof of 
Proposition~\ref{prop-3.1} for the case
$p\geq 1$.
The case
$p=0$
can be done by  similar methods.

A general outline of what we are going to do in this section
is the following. We start with a neighborhood
$U$
of the disk
$D^p$
in
$\MKN$, 
then we write the function
$f$
in a normal form on a smaller neighborhood
$U^\prime\subset U$ 
and then we construct a deformation
$\{ f_\tau\}$
of
$f$
which does not change 
$f$
outside of
$U^\prime$.

\bigskip
\noindent
{\bf Notational convention.}

\smallskip
\noindent
From now on let
$t=(t_1,\ldots ,t_s)$,
$y=(y_1,\ldots ,y_p)$,
$z=(z_1,\ldots ,z_{n-s-p})$,
$\xi=(\xi_1,\ldots ,\xi_N)$,
${\bar\xi}=(\xi_2,\ldots ,\xi_N)$.
As before by
$B_r^l$
we denote an
$l{\hbox{\rm -dimensional}}$
ball of radius
$r$
with the center at 
$0\in {\bf R}^l$
and we omit the lower index for unit balls. 
\bigskip

\subsection{Auxiliary constructions}
\label{pf-prop-3.1-aux}

\begin{lem}
\label{lem-4.1} 
For an arbitrarily small neighborhood 
$U$ 
of the disk 
$D^p$
in 
$M\times\KN$
there exists a diffeomorphism 
$j: {\R}^{n+N}\to U$ 
such that the following conditions are fulfilled:

\smallskip
\smallskip
\noindent
(A) 
on a neighborhood
${U^\prime}\subset U$  
of 
$j^{-1}(D^p)\subset 
\RnN=\Rp\times\Rnsp\times {\bf R}^s\times\R\times {\bf R}^{N-1}$
one has

{\bigskip
\noindent
$\displaystyle f\circ j (y,z,t,\xi_1,{\bar\xi}) =$
\nopagebreak
\smallskip
\noindent
$\displaystyle \lambda_{s+3}(y,z)\{ \xi_1^{s+3} +\theta (y,z)\xi_1^{s+1}
+ \sum_{i=1}^s t_i\xi_1^i\} + \lambda_0 (y,z,t) -
\sum_{i=2}^{2+k}\xi_i^2 +\sum_{i=3+k}^N\xi_i^2$
\bigskip
}

\noindent
for an integer number
$k$,
$0\leq k\leq N-1$,
and some functions
$\lambda_0(y,z,t)$,
$\lambda_{s+3}(y,z)$,
$\lambda_{s+3}(y,z)>0$
on
$U^\prime$,
$\theta(y,z) ={\theta^\prime}(y,z)/\lambda_{s+3}(y,z)$,
where a function
$\theta^\prime (y,z)$
satisfies the following conditions:

\smallskip
(A1)
$\theta^\prime$
is positive on
$({\R}^p \setminus B^p)\times B^{n-s-p}$,
negative on
$B^p \times (\Rnsp\setminus B^{n-s-p})$,

(A2)
$\theta^\prime$
has a nondegenerate critical point of index
$p$
at zero and it is the only critical
point of
$\theta^\prime$;
\smallskip

\smallskip
\noindent
(B1)
$j$
maps the fibers of the projection
${\R}^{n+N}\to {\R}^n$
into the fibers of the projection
$\MKN\to M$;

\smallskip
\noindent
(B2)
$j$
maps the line field parallel to the axis
$\xi_1$ 
on 
$\RnN$
into the line field
$l =  ker\ dP$,
where 
$P$
is the restriction of the projection
$\MKN\to M$
on
$V_L$;

\smallskip
\noindent
(B3)
$j^{-1}(D^p)\subset {\R}^p\times 0\times 0\times 0\times 0$;

\smallskip
\noindent
(B4)
$j^{-1}(V)\subset \Rp\times\Rnsp\times{\bf R}^s\times {\bf R}
\times 0$;

\smallskip
\noindent
(B5)
$j^{-1}(D^p\cup V_{s+1}), j^{-1}(D^p\times[-1,1])\subset \Rp\times\Rnsp
\times 0\times 0\times 0\subset j^{-1}(V_{s-1})$;

\smallskip
\noindent
(B6)
$j^{-1}(V_{s+1})\subset \{(y,z,0,0,0)\ | \  (y,z,0,0,0)\in 
\Rp\times\Rnsp\times 0\times 0\times 0\subset\RnN, 
{\theta^\prime} (y,z)=0\}$;

\smallskip
\noindent
(B7)
$j^{-1}(V_s)\cap (\Rp\times\Rnsp\times 0\times 0\times 0) = 
V_{s+1}$;

\smallskip
\noindent
(B8) 
$j^{-1}(V_{s-1})\subset 
{\bf R}^p\times {\bf R}^{n-s-p}\times 0\times \R\times {\bf R}^{N-1}$.
\end{lem}

\noindent
We prove Lemma~\ref{lem-4.1} in Section~\ref{pf-lemm-4.1}.

\begin{rem}
\label{rem-4.1}
{\rm
Essentially Lemma~\ref{lem-4.1} is about a 
{\it parametric Morin normal form} 
(\cite{M}, see also \cite{AGV}) for singularities of
the type
${\Sigma}^{11\ldots 1}.$
Similar result was implicitly used in the 
proof of Proposition 3.8 in 
\cite{El}.
Another result of this kind can be found in
\cite{Ig}.
}
\end{rem}
\bigskip

Now we can write explicitly the function
$f$
in a small neighborhood
$U^\prime$
of
$D^p$
and our next step is to find a necessary
$C^1{\hbox{\rm -small}}$
deformation of 
$f$
which changes
$f$
{\it 
only inside
}
$U^\prime$.
For this purpose we need the following construction
which takes the rest of this section.

Let us use Lemma~\ref{lem-4.1} 
and let us identify
$U$
with
$\RnN$
(so from now on we omit
$j$
everywhere).
We can assume that 
${U^\prime}=B^p\times B^{n-s-p}\times B^s_{\varepsilon_1}\times
B^1_{\varepsilon_2}\times B^{N-1}$
(for some
$\varepsilon_1, \varepsilon_2$).
We can also assume that the function
$\theta\leq C$
on
$U^\prime$
for some positive constant
$C$.

For any 
$\epsilon$,
$0<\epsilon<\varepsilon_2$,
there exist a number
$\delta(\epsilon)>0$
and an open ball
$B_{r(\epsilon)}^s$,
$r(\epsilon)<\varepsilon_1$,
such that for any
$a$,
$-\delta\leq a\leq C$,
and any
$t\in B_{r(\epsilon)}^s$
all critical points of the function
\[x\to t_1 x^1+\ldots+t_s x^s+ax^{s+1}+x^{s+3}\]
are contained in
$B_{\epsilon/2}^1$.
Let
$B_R^{n-s}\subset B^p\times B^{n-s-p}$
be a ball containing
$D^p\subset \Rp\times\Rnsp\times 0\times 0\times 0$
such that on
$B_R^{n-s}$
one has
$\theta (y,z)\geq -\delta(\epsilon)$.
Then for any 
$(y,z,t)\in B_R^{n-s}\times B_{r(\epsilon)}^s$
all the critical points of the function
\[\xi_1\to \xi_1^{s+3} +{\theta} (y,z)\xi_1^{s+1}
+ \sum_{i=1}^s t_i\xi_1^i\]
lie inside
$B_{\epsilon/2}^1$.
Let us choose an open ball
$B_{R_1}^{n-s}$,
$D^p\subset B_{R_1}^{n-s}\subset B_R^{n-s}$, 
and let us define a  family of functions
$\{ \theta_\tau \}:{\bf R}^p \times 
{\bf R}^{n-s-p} \to \R$,
$0\leq\tau\leq 1$,
with the following properties:

\medskip
(I) 
$\theta_0=\theta$;

\smallskip
(II) for any
$\tau$,
$0\leq\tau\leq 1$,
$\theta_\tau =\theta$
outside
$B_{R_1}^{n-s}$;

\smallskip
(III) for any
$\tau$,
$0\leq\tau\leq 1$,
$\theta_\tau (y,z)\geq -\delta(\epsilon)$
on
$B_R^{n-s}$;
 
\smallskip
(IV) the function
$\theta_1$
is negative on
$B_R^{n-s}$.
\medskip

\smallskip
Hence for any
$\tau$
and any 
$(y,z,t)\in B_R^{n-s}\times B_{r(\epsilon)}^s$
all the critical points of the function
\[\xi_1\to \xi_1^{s+3} +\theta_\tau (y,z)\xi_1^{s+1}
+ \sum_{i=1}^s t_i\xi_1^i\]
lie inside
$B_{\epsilon/2}^1$.
Then we can easily
find a family of functions
$\{ g_\tau\}:{U^\prime}\to \R$,
$0\leq\tau\leq 1$,
such that
$g_0 = f$
on
$U^\prime$
and for any
$\tau$

\smallskip
\noindent
(i) 
$g_\tau = f$
outside 
$B_{R_1}^{n-s}\times B_{r(\epsilon)}^s\times B_\epsilon^1$;

\smallskip
\noindent
(ii)  
\bigskip
\noindent
$g_\tau (y,z,t,\xi_1,{\bar\xi}) = \lambda_{s+3}(y,z)(\xi_1^{s+3} +
\theta_\tau (y,z)\xi_1^{s+1}
+\sum_{i=1}^s t_i\xi_1^i) + \lambda_0 (y,z,t) -
\sum_{i=2}^{2+k}\xi_i^2 +\sum_{i=3+k}^N\xi_i^2$
on 
$B_R^{n-s}\times B_{r(\epsilon)/2}^s\times B_{\epsilon/2}^1$;

\smallskip
\noindent
(iii) 
${\partial g_\tau}/{\partial \xi_1}\neq 0$
outside 
$B_R^{n-s}\times B_{r(\epsilon)/2}^s\times B_{\epsilon/2}^1$.
\bigskip

\subsection{End of the proof of Proposition~\ref{prop-3.1}}
\label{pf-prop-3.1-end}

Let us define
$f_\tau$
as
$f$
outside
$U$
and as
$g_\tau$
inside 
$U$.
Taking a sufficiently small
$\epsilon$
we can make the deformation
$\{ f_\tau \}$ 
arbitrarily
$C^1{\hbox{\rm -small}}$
on
$U$.

Let us check  the condition (iii) from Proposition~\ref{prop-3.1}. It follows from Lemma~\ref{lem-4.1} 
and from the construction of
$\{ g_\tau \}$
that we only need to check that the 
$\alpha{\hbox{\rm -chain}}$ 
generated by the function
$\Psi_1 :\Rn\times\R\to \R$
is obtained from the 
$\alpha{\hbox{\rm -chain}}$ 
generated by the function
$\Psi_0:\Rn\times\R\to \R$
by a direct surgery along the disk
$D^p$,
where
\[\Psi_\tau (y,z,t,\xi_1) = \sum_{i=1}^{s-1}t_i\xi_1^i + 
\theta_\tau(y,z)\xi_1^{s+1} + \xi_1^{s+3} \ (0\leq\tau\leq 1).\]
In that case we indeed have the required surgery
and therefore the condition (iii) holds.

Now let us check the condition (ii). Lemma~\ref{lem-4.1} 
and the construction of 
$\{ g_\tau \}$ 
imply that we only need to check that the function 
$\Psi_\tau :\Rn\times\R\to\R $
generates an embedded Lagrangian submanifold in
${T^\ast}\Rn$.
For any fixed 
$(y,z,t)$
and any
$x_j\in\R\ (j=1,2)$
such that 
${\partial\Psi_\tau}/{\partial x} (y,z,t,x_j)=0$
we have that
\[{{\partial\Psi_\tau }\over {\partial t_1 }} (y,z,t,x_j) 
= x_j \ (j=1,2).\]
Then necessarily 
$x_1=x_2$.
Hence the Lagrangian submanifold 
$L_\tau^\prime\subset T^\ast L$ 
generated by 
$\left. f_\tau\right|_U$
is embedded.
Therefore, since  we change
$f$
by an arbitrarily 
$C^1{\hbox{\rm -small}}$
deformation in an arbitrarily small neighborhood of
$D^p$,
we get that the whole Lagrangian submanifold
$L_\tau$
generated by
$f_\tau$
is embedded. 
Thus
$L_\tau$
is  Lagrangian isotopic to
$L$
{\it 
in the class of embedded
Lagrangian submanifolds.
}
Since this isotopy changes
$L$
only in a small neighborhood of the disk
$D^p$,
it does not change the periods of the canonical
$1{\hbox{\rm -form}}$
$pdq$
on
$L$.
Hence
$L_\tau$
is
{\it
Hamiltonian
isotopic}
to
$L$.

Let us denote the Hamiltonian flow which deforms 
$L$
into
$L_1$
by
$\{ h_\tau\}$.
Let us fix a point 
$x_0\in L$
and let us consider the family of the action functions
$\{A_\tau\}$, 
defined by the formula
$A_\tau (x) = \int_{\gamma_\tau} pdq :L_\tau\to {\bf R}$,
where
$\gamma_\tau$
is a path in
$L_\tau$
connecting
$h_\tau (x_0)$
with
$x$.
The function
$A_\tau$
is correctly defined because the Lagrangian submanifold
$L_\tau$
is exact.
Up to a global constant (independent of
$\tau$), 
$A_\tau$
is equal to
$f_\tau$
for any
$\tau$.
Hence 
${\mid\int_{h_\tau(\gamma)} pdq \mid}$
can be made uniformly bounded by an arbitrarily small number
independent of
$\tau$
and a path 
$\gamma$
in
$L$.
This finishes the proof of 
Proposition~\ref{prop-3.1}.
\b

\subsection{Proof of Lemma~\ref{lem-4.1}}
\label{pf-lemm-4.1}

Let us pick an arbitrary function
$\theta^\prime (y,z)$
that satisfies conditions (A1),(A2) of 
Lemma~\ref{lem-4.1}.
Since 
$D^p$
is a part of a basis of the direct surgery and since it projects 
into an 
{\sl
embedded
}
disk in
$M$,
one finds that for an arbitrarily small neighborhood
$U$
(diffeomorphic to
$\RnN$)
of 
$D^p$
in
$\MKN$
there exists a diffeomorphism
$j:\RnN\to U$
which satisfies conditions (B1)--(B8) of Lemma~\ref{lem-4.1}.
We call such a pair
$(U,j)$
a
{\it
good chart.}
 
We call a (local)  change of coordinates in
a neighborhood of
$j^{-1}(D^p)\subset
\RnN=\Rp\times\Rnsp\times {\bf R}^s\times\R\times {\bf R}^{N-1}$
{\it good,}
if it has the form
\[(y,z,t,\xi_1 ,{\bar\xi})\to (y,z,t,
S_1(y,z,t,\xi_1,{\bar\xi}),S_2(y,z,t,\xi_1,{\bar\xi})),\]
$S_1:\RnN\to\R$, 
$S_1(y,z,t,0,{\bar\xi})=0$,
$S_2:\RnN\to {\bf R}^{N-1}$, 
$S_2(y,z,t,\xi_1,0)
\break
= 0\in {\bf R}^{N-1}.$
Obviously, a good change of coordinates maps a good chart into
another good chart.

We break  the proof of Lemma~\ref{lem-4.1} into several steps (Lemmas~\ref{lemA}--\ref{lemE}).
Lemma~\ref{lem-4.1}  
follows immediately from Lemma~\ref{lemE}.

\begin{lem}
\label{lemA}
There exists a good chart
$(U,j)$
such that the function
$f\circ j:\RnN\to\R$
has the form
\[f\circ j(y,z,t,\xi )=\varphi (y,z,t,\xi_1)-\sum_{i=2}^{2+k}\xi_i^2 +
\sum_{i=3+k}^N \xi_i^2\]
in a neighborhood of
$j^{-1}(D^p)$
for some function
$\varphi$
and some integer  
$k$,
$0\leq k\leq N-1.$
\end{lem}

\begin{lem}
\label{lemB}
The function
$\varphi$
in 
Lemma~\ref{lemA} 
can be reduced to the form
\setcounter{equation}{0}
\renewcommand{\theequation}{\fnsymbol{equation}}
\begin{equation}
\varphi (y,z,t,\xi_1)=
\sum_{i=0}^{s+3}\lambda_i (y,z,t)\xi_1^i,
\end{equation}
where
$\lambda_{s+3}(y,z,t)\neq 0$
for any
$(y,z,t)$,
and if
$(y,0,0,0,0)\in j^{-1}(\partial D^p)$
then
$\lambda_i (y,0,0)=0$,
$i=1,\ldots ,s+2.$
\end{lem}

\begin{lem}
\label{lemC}
One can assume that in 
Lemma~\ref{lemB} 
\[\lambda_{s+3}(y,z,t)\equiv \pm 1,\ 
\lambda_{s+2}(y,z,t)\equiv 0,\ 
\lambda_i (y,z,0)=0,\  i=1,\ldots ,s,\]

and
\[\lambda_{s+1} (y,z,0)=0 \ {\it iff}\  (y,z,0,0,0)\in j^{-1}(\partial D^p).\]
\end{lem}

\begin{lem}
\label{lemD}
One can assume that in Lemma~\ref{lemA}
\[\varphi (y,z,t,\xi_1)=\pm\xi_1^{s+3}+\lambda_{s+1} (y,z)\xi_1^{s+1} + 
t_s\xi_1^s +\ldots +t_1\xi_1 +\lambda_0 (y,z,t),\]
where
$\lambda_{s+1} (y,z)=0 \ {\it iff}\ (y,z,0,0,0)\in j^{-1}(\partial D^p).$
\end{lem}

\begin{lem}
\label{lemE}
There exists a good chart
$(U,j)$
such that

\bigskip
\noindent
$\displaystyle f(y,z,t,\xi_1,{\bar\xi}) =$

\smallskip
\noindent
$\displaystyle \pm\lambda_{s+3}\xi_1^{s+3} +{\theta^\prime} (y,z)\xi_1^{s+1} + 
t_s\xi_1^s +\ldots
+t_1\xi_1 +\lambda_0 (y,z,t) - 
\sum_{i=2}^{2+k}\xi_i^2 +\sum_{i=3+k}^N\xi_i^2$
\bigskip

\noindent
in a neighborhood of
$D^p$.
\end{lem}
\bigskip

In the proofs of the lemmas we use parametric versions
of ``the Morse lemma with parameters'' (\cite{Ar1}, see also \cite{AGV}) and
of the Tougeron theorem (\cite{T}, see also \cite{AGV}). The
parametric versions  that we need for our case (when the corank of
the function is equal to 1 for any value of the parameter) 
can be easily proved along the lines of the
proofs of the standard ``nonparametric'' versions (see \cite{AGV}).
For example, the exact statement of the version of the Tougeron
theorem that we need is as follows.

\bigskip
\noindent
{\bf Parametric Tougeron theorem.}
{\it
Suppose we have two families
$\{f_a\}:(\R,0)\to (\R,0)$,
$\{g_a\}:(\R,0)\to (\R,0)$,
of germs of smooth functions, where the parameter 
$a$ 
runs over a compact manifold. Suppose that for any
$a$
the singularities of
$f_a$
and
$g_a$
at zero have the same Milnor number
$\mu_a$
and the same
$(\mu_a+1){\hbox{\it -jet}}$.
Then there exists a smooth family 
$\{S_a\}:(\R,0)\to(\R,0)$
of germs of diffeomorphisms of 
$\R$
which takes the family
$\{ f_a\}$
to the family
$\{ g_a\}$.
}
\bigskip

From now on, given a good chart
$(U,j)$,
we identify 
$\RnN$
with
$U=j(\RnN).$

\bigskip
\noindent
{\bf Proof of Lemma~\ref{lemA}.}
For any point
${\bf x}=(y,z,t,\xi_1, 0)\in V$
the corank of the second differential of the function
$f(y,z,t,\cdot ,\cdot):\R\times {\bf R}^{N-1}\to\R$
at
${\bf x}$
is equal to 1 (because 
we do not allow any
$\Sigma^2{\hbox{\rm -singularities}}$
of tangency) 
and
the kernel of its second differential is given by the
axis
$\xi_1$
in
$\R\times {\bf R}^{N-1}.$
Therefore 
$|{\partial}^2 f/{\partial{\xi_{i_1}}}{ \partial{\xi_{i_2}}}({\bf x})|\neq 0$,
$i_1,i_2=2,\ldots ,N$,
and thus the signature of the second differential of the function
$f(y,z,t,\cdot ,\cdot):\R\times {\bf R}^{N-1}\to\R$
is constant in a neighborhood of the disk
$D^p\subset\RnN.$
By the parametric version of the ``Morse lemma with
parameters'' there exists a good change of coordinates of
the type  
$(y,z,t,\xi_1,{\bar\xi})\to (y,z,t,\xi_1, S_2(y,z,t,\xi_1,{\bar\xi}))$
which reduces
$f$
to the necessary form. 
\b
\bigskip

\bigskip
\noindent
{\bf Proof of Lemma~\ref{lemB}.}
At any point
${\bf x}=(y,z,t,\xi_1,0)\in V_k\setminus V_{k+1}$,
$0\leq k\leq s+1$,
we have that
${{\partial}^i f}/{{\partial\xi_1^i}} ({\bf x})=0$,
$i=1,\ldots ,k+1$,
and
${{\partial}^{k+2} f}/{\partial\xi_1^{k+2}} ({\bf x})\neq 0.$
A sufficiently small neighborhood of the disk
$D^p$
intersects 
$V_i$
only if
$i\leq s+1$.
Therefore for any point
$(y,z,t,0,{\bar\xi})$
close enough to
$D^p$
the Milnor number of the function
$f(y,z,t,\cdot,{\bar\xi}):\R\to\R$
at zero is less or equal to
$s+2$.
Hence the parametric version of the Tougeron theorem provides
us with a good change of coordinates of the type
$(y,z,t,\xi_1 ,{\bar\xi})\to (y,z,t,
S_1(y,z,t,\xi_1,{\bar\xi}),{\bar\xi})$,
which reduces
$\varphi$
to the necessary form.

For any
$(y,0,0,0,0)\in \partial D^p$
we have
$\lambda_i(y,0,0)=0$,
$i=1,\ldots ,s+2$,
$\lambda_{s+3}(y,0,0)\neq 0$.
Moreover, the function
$\lambda_{s+3}$
does not change sign along 
$\partial D^p$
(for
$p>1$
this is obvious;
for
$p=1$
it follows from Definition~\ref{def-1.7}
and from 
Remarks~\ref{rem-1.2a},\ref{rem-1.3}).
Therefore there exists a function
$\tilde\lambda_{s+3}$
on a neighborhood of
$D^p$
which does not vanish on that neighborhood and which 
coincides with
$\lambda_{s+3}$
on a neighborhood of
$\partial D^p$.
Applying the parametric version of the Tougeron theorem once again
by a good change of coordinates we can reduce 
$\varphi$
to the form
$(\ast)$
but with the function
$\tilde\lambda_{s+3}$
in the place of
$\lambda_{s+3}.$
\b
\bigskip

\bigskip
\noindent
{\bf Proof of Lemma~\ref{lemC}.}
By a shift and a dilation of 
$\xi_1$
we obtain the first two conditions.
The third condition follows from the fact that
any point
$(y,z,0,0,0)$
belongs to
$V_{s-1}$.
The last condition holds because the disk
$D^p$
lies inside
$V_{s-1}$
and intersects
$V_{s+1}$
normally.
\b
\bigskip

\bigskip
\noindent
{\bf Proof of Lemma~\ref{lemD}.}
According to 
Lemma~\ref{lemC}, 
we can assume that we already have
a good chart for which
$\varphi (y,z,t,\xi_1) = \varphi_1 (y,z,t,\xi_1) +\lambda_0 (y,z,t),$
where
\[\varphi_1 (y,z,t,\xi_1)=\pm\xi_1^{s+3} + 
\sum_{i=1}^{s+1} \lambda_i (y,z,t)\xi_1^i.\]
We also can assume that the original generating function
$f$
is generic (otherwise we could have 
${C^\infty}{\hbox{\rm -slightly}}$ 
perturbed it at the very beginning). Therefore
we can  assume that for any
${\bf x}\in D^p$
the mapping
\[(y,z,t,\xi_1)\to (y,z,t,\varphi_1 (y,z,t,\xi_1)):{\bf R}^{s+1}\to 
{\bf R}^{s+1}\]
is 
$RL{\hbox{\rm -stable}}$
at 
${\bf x}$
(see \cite{AGV}). In this case,
since according to Lemma~\ref{lemC},
\[\partial\lambda_i/\partial y_j (y,z,0)=0,
\partial\lambda_i/\partial z_k (y,z,0)=0,  
\]
\[
i=1,\ldots ,s,
j=1,\ldots ,p, k=1,\ldots ,n-s-p,\]
the theorem from \cite{AGV}, ch. 9.5, implies that        
for the points
$(y,z,t,0,0)$
close enough to
$D^p$
\[{\rm det}\ \| \partial\lambda_i (y,z,t)/\partial t_j\|\neq 0\
(i,j = 1,\ldots ,s-1).\]
So we can make a change of coordinates
\[(y,z,t_1,\ldots,t_{s-1},t_s,\xi_1,{\bar\xi})\to
(y,z,\lambda_1 (y,z,t),\ldots ,\lambda_{s-1} (y,z,t),
t_s, \xi_1,{\bar\xi})\]
in a neighborhood of
$D^p$
so that the new chart is good and in the new coordinates
\[\varphi_1 (y,z,t,\xi_1) = \xi_1^{s+3} + \lambda_{s+1}(y,z,t)\xi_1^{s+1}
+\lambda_s (y,z,t)\xi_1^s +\sum_{i=1}^{s-1}t_i\xi_1^i.\]
Applying the parametric version of the Tougeron theorem and using
the fact that a good chart satisfies condition (B8) of Lemma~\ref{lem-4.1}
we can reduce
$\varphi_1$
by a good change of coordinates to the following form:
\[\varphi_1 (y,z,t,\xi_1) = \xi_1^{s+3} + \lambda_{s+1}(y,z)\xi_1^{s+1}
+\lambda_s (y,z,t_s)\xi_1^s +\sum_{i=1}^{s-1}t_i\xi_1^i,\]
where
$\lambda_{s+1} (y,z)=0$ 
if and only if
$(y,z,0,0,0)\in V_{s+1}$
and
$\lambda_s (y,z,t_s)=0$
if and only if
$t_s=0$. 

Since 
$\lambda_s (y,z,t_s)=0$
if and only if
$t_s=0$
the function
$\partial\lambda_s/\partial t_s$
does not change sign along
$D^p$.
Therefore after a 
$C^\infty{\hbox{\rm -small}}$
perturbation of
$\lambda_s$
(and thus a
$C^\infty{\hbox{\rm -small}}$
perturbation of the generating function
$f$)
in a small neighborhood of
$D^p$ 
we can achieve that
$\partial\lambda_s/\partial t_s\neq 0$
along
$D^p$. 
Then we can make the change of coordinates
\[(y,z,t_1,\ldots,t_{s-1},t_s,\xi_1,{\bar\xi})\to
(y,z,t_1,\ldots ,t_{s-1},\lambda_s (y,z,t),\xi_1,{\bar\xi})\]
in a neighborhood of 
$D^p$
which would reduce 
$\varphi$
to the necessary form.
\b
\bigskip

\bigskip
\noindent
{\bf Proof of Lemma~\ref{lemE}.}
One only needs to use Lemma~\ref{lemD} and to observe that 

\noindent
1) according to our construction and the previous lemmas, 
the zero sets of the functions
$\theta^\prime$
and
$\lambda_{s+1}$
coincide;

\noindent
2) according to the definition of a basis of a direct
surgery, the sign of the function
$\pm\lambda_{s+3}$
must be the same as the sign of
the function
$\lambda_{s+1}$
inside
$D^p$.
\b

\begin{figure}
\centerline{\psfig{figure=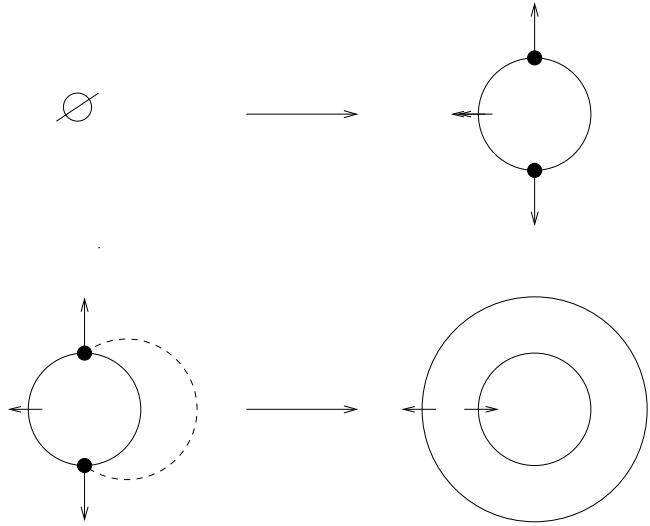,height=80mm}}
\caption{Creation of a double fold by direct surgeries.}
\label{fig10}
\end{figure}

{\section{Double folds: creation and normal forms}
\label{double-folds-pfs}
\nopagebreak
\bigskip
\noindent
{\bf A proof of Propositions~\ref{prop-2.2} and \ref{prop-2.2a}.}
\nopagebreak
Let 
$x$
be a point inside the ball
$B$.
We create a double fold near
$x$
by applying two subsequent direct surgeries and using Proposition~\ref{prop-3.1}.
The surgeries for one of the two possible double folds
in the two-dimensional case are shown on Fig.~\ref{fig10}.
}

More precisely, we first apply a direct surgery of order 1 and index 0 inside
a neighborhood of
$x$
in
$B$
and create a sphere of folds 
$S_1^{n-1}$
with a sphere of cusps
$S_2^{n-2}\subset S_1^{n-1}\subset B$
on it.
Let
$B_1\subset L$
be the ball bounded by
$S_1^{n-1}$.
The vector field 
$v_2$
on
$S_2^{n-2}$
necessarily points outside of 
$B_1$
and either one of two possible  Maslov coorientations of
$S_1^{n-1}$
in
$L$
can be realized.
Then we consider an embedded disk
$D^{n-1}\subset L$
which lies in a small neighborhood of the ball 
$B_1$
so that
$\partial D^{n-1}= S_2^{n-2}$
and the interior of
$D^{n-1}$
lies outside of the ball
$B_1$.
We apply a direct surgery of index 
$n-1$
and order
$1$
along
$D^{n-1}$
and create a double fold that we need.
According to 
Proposition~\ref{prop-3.1}, 
the direct surgeries above can be
realized by an appropriate deformation of
$L$.
\b
\bigskip

Alternatively, we can present an explicit construction of a double fold
using a deformation of generating functions.

\bigskip
\noindent
{\bf Another proof of Propositions~\ref{prop-2.2} and \ref{prop-2.2a}.}
We present a direct proof of 
Proposition~\ref{prop-2.2}. One can prove 
Proposition~\ref{prop-2.2a} using the same construction.

Again, let
$x$
be a point inside the ball 
$B$.
Using the Darboux theorem for Lagrangian vector bundles (see e.g. \cite{AGV})
we find some symplectic coordinates near
$x$
in
$T^\ast L$
so that near
$x$
the submanifold
$L$
is identified with a neighborhood of the origin in the zero-section of
$T^\ast {\bf R}^n$
and 
${\cal F}$
is identified with the foliation of
$T^\ast {\bf R}^n$
by the leaves of the canonical Lagrangian projection
$T^\ast {\bf R}^n \to {\bf R}^n$.
Thus we can assume that
$B=B^n\subset {\bf R}^n$
is the standard unit ball and then
it is enough for us to construct a Hamiltonian isotopy of the ball 
$B$
in
$T^\ast {\bf R}^n$
with the required properties
(the singularities of tangency to
$\cal F$
are now viewed as singularities of the Lagrangian projection
$T^\ast {\bf R}^n \to {\bf R}^n$).
We can construct this Hamiltonian isotopy of
$L$
by constructing a deformation of generating functions.

We start with a function on
${\bf R}^n\times {\bf R}$
which generates the zero-section
$L\subset {\bf R}^4$.
We then construct a deformation of this generating function
$f_0$
which in a neighborhood of the origin deforms the flat
$V_L = {\bf R}^n\times 0 \subset {\bf R}^n\times {\bf R}$
into a ``mushroom'' with respect to the ``vertical'' fibers of the projection
${\bf R}^n\times {\bf R}\to {\bf R}^n$ as on Fig.~\ref{fig11}.
The construction proceeds as follows.

Given an
$\epsilon\in (0,1)$
and
$r_1, r_2,r_3,r_4$,
$0<r_1<r_2<r_3<r_4<1$,
let
$\{ t_a(r)\}:\R\to\R$,
$-\epsilon\leq a \leq 0$,
be a smooth family of functions which satisfies the following 
properties for any
$a$:

\medskip
\noindent
(i) 
$t_a(r)\equiv 0$
outside of the interval
$[r_1,r_4]$;

\smallskip
\noindent
(ii)
$t_a(r)$
is a linear function on
$[r_2, r_3]$
and
$t_a (r_2) =  100a$,
$t_a (r_3) =  -100a$.
\medskip

\smallskip
For sufficiently small
$\epsilon, r_1, r_2, r_3, r_4$
we can construct a family of functions
$\{ f_a\}:B^n\times B^1\to {\bf R}$,
$-\epsilon\leq a\leq 0$,
some  balls
$B^n_{r_5}$,
$B^n_{r_6}$,
$r_4<r_5<r_6<1$,
and
$B^1_{\delta_1}$,
$B^1_{\delta_2}$,
$0<\delta_1<\delta_2<1$,
so that for any 
$a$,
$-\epsilon\leq a\leq 0$,
the following conditions hold:

\medskip
\noindent
(i) inside
$B^n_{r_5}\times B^1_{\delta_1}$
the function
$f_a$
is given by the formula
\[f_a (y,\xi) = \xi^4 +a\xi^2 +t_a(r)\xi ,\]
$y=(y_1,\ldots ,y_n)\in B^n$,
$r=r(b)=y_1^2+\ldots+y_n^2$,
$\xi\in B^1$;

\smallskip
\noindent
(ii) 
$f_a (y,\xi) = \xi^2  ,$ 
outside 
$B^n_{r_6}\times B^1_{\delta_2}$;

\smallskip
\noindent
(iii) for any
$a$
all the critical points of
$f_a$
with respect to the variable
$\xi$
lie inside
$B^n_{r_5}\times B^1_{\delta_1}$.
\medskip

\smallskip
Consider a smooth family of functions
$\{ g_\tau = f_{-\epsilon\tau} \}$,
$0\leq\tau\leq 1$,
on
$B^n\times {\bf R}$.
For any
$\tau$
the function
$g_\tau$
generates an
{\it embedded}
Lagrangian submanifold in
$T^\ast {\bf R}^n$
which, if we started with a sufficiently small
$\epsilon$,
coincides with 
$L$
outside of an arbitrarily small neighborhood of
$0$
in
$T^\ast {\bf R}^n$.
One easily checks that
$g_1$
generates a Lagrangian submanifold
$L_1\subset T^\ast {\bf R}^n$
such that the only singularities of the Lagrangian projection on
$L_1$
are some two small spheres
$S_1, S_2\subset L_1$
of folds with the opposite Maslov coorientations.
If we change
$t_a(r)$
to
$-t_a(r)$
in our construction, the Maslov coorientations of
$S_1$ 
and
$S_2$
change simultaneously and thus we can always get a necessary
combination of them.

It is not hard to see that the deformation
$\{ g_\tau\}$ 
gives us  a Lagrangian isotopy
$\{ L_\tau\}$
of
$L$
in
$T^\ast {\bf R}^n$
which is identical outside
$B^n$.
This Lagrangian isotopy is actually a Hamiltonian  isotopy 
$\{L_\tau = h_\tau (L)\}$
because 
under our deformation the periods of the canonical
$1{\hbox{\rm -form}}$
$pdq$
on
$L$
do not change.

By an appropriate choice of
$\epsilon, r_1, r_2, r_3, r_4$
the deformation
$\{ g_\tau\}$ 
can be made arbitrarily
$C^1{\hbox{\rm -small}}$.
Therefore the Hamiltonian isotopy
$\{ h_\tau\}$
can be chosen in such a way that it is arbitrarily
$C^0{\hbox{\rm -small}}$
and so that 
$\mid\int_{h_\tau(\gamma)} pdq \mid$
is uniformly bounded by an arbitrarily small number independent of
$\tau$
and a path 
$\gamma$
in
$L$.
\b
\bigskip

Now let us the discuss normal forms for double folds mentioned in 
Remark~\ref{d-fold-n-form}. 

Suppose that the zero-section
$L\subset T^\ast L$
(resp.
$L^\prime\subset T^\ast L^\prime$)
has a double fold along some spheres 
$S_1$, 
$S_2$
(resp.
$S_1^\prime$,
$S_2^\prime$)
with respect to a Lagrangian foliation
$\cal F$
(resp.
${\cal F}^\prime$.
Suppose that the Maslov coorientations of
$S_1$
and
$S_2$
form the same combination as for
$S_1^\prime$
and
$S_2^\prime$.
Let
$B$
(resp.
$B^\prime$)
be a ball in 
$L$
(resp.
$L^\prime$)
bounded by 
$S_2$
(resp.
$S_2^\prime$),
the larger of the two spheres of folds.

\begin{prop}
\label{prop-dfnf}
There exists a symplectomorphism 
$F$
between a neighborhood of
$B$
in
$T^\ast L$
and a neighborhood of
$B^\prime$
in
$T^\ast L^\prime$
which  maps the zero-section
$L\subset T^\ast L$
into the zero-section
$L^\prime \subset T^\ast L^\prime$
and
$\cal F$
into
${\cal F}^\prime$
(so it also has to map
$S_1$
into
$S_1^\prime$
and
$S_2$
into
$S_2^\prime$).

In particular, for any 
$L$
with a double fold along 
$S_1$
and
$S_2$
we can find such a symplectomorphism
$F$
which identifies a small neighborhood of the ball 
$B$
with a small neighborhood of one of the explicit double folds
constructed in the second proof of 
Proposition~\ref{prop-2.2} above.

\end{prop}

\bigskip
\noindent
{\bf Proof of 
Proposition~\ref{prop-dfnf}.}
First, let us show that we can construct the symplectomorphism
$F$
in a neighborhood of each sphere of folds
$S_i$.

Let 
$S$ 
be one of the spheres
$S_1$
and
$S_2$.
Consider a small neighborhood 
$\widetilde U$
of the sphere 
$S$
in
$T^\ast L$.
We can choose an almost complex structure
$J$
on
$\widetilde U$
so that
$S$
lies in a single leaf 
$A$
of the Lagrangian foliation
$J({\cal F})$
of
$\widetilde U$.
A small neighborhood of
$S$
in
$A$
(still denoted by
$A$)
is an 
$n{\hbox{\rm -dimensional}}$
Lagrangian spherical annulus.
Hence, using the Lagrangian neighborhood theorem 
\cite{We1},\cite{W}
we can symplectically identify a small neighborhood of
$A$
in
$T^\ast L$
with a small neighborhood 
(still denoted by 
$\widetilde U$)
of the zero-section
$A$
in
$T^\ast A$.
Let us identify
$A$
with an 
$n{\hbox{\rm -dimensional}}$
spherical annulus in
${\bf R}^n\subset {\bf R}^{2n}$
(still denoted by
$A$)
given by the equations
$p_1 =p_2=\ldots = p_n =0, 2-\delta < \|q\| < 2+\delta$,
where
$\delta$
is very small.
According to 
\cite{G},
we can find a generating function
$g$
for the Lagrangian submanifold
$L\cap {\widetilde U}$.
This function 
$g$
will be defined on an open subset of  
$A\times {\bf R}^N$
for some 
$N$.
By the same argument as in the proof of
Lemma~\ref{lem-4.1}
(see Section~\ref{pf-lemm-4.1})
we can reduce 
$g$
to a normal form
$g_{\alpha ,\beta ,Q} (q,\xi) = \xi_1^3 + 
\alpha ({\bf r})\xi_1 + \beta (q) + Q(\xi_2, ...,\xi_n)$.
Here
$q$
runs over a small neighborhood of
$S$
in
$A$,
$\xi = (\xi_1,\ldots, \xi_N)$
are the coordinates on a neighborhood of zero in
${\bf R}^N$,
${\bf r} = q_1^2+\ldots +q_n^2$
and
$Q$
is a nondegenerate quadratic form.
The function
$\alpha ({\bf r})$
has a simple zero along
$S$
and the sign of
$\partial \alpha/\partial {\bf r}$
on
$S$
depends on the Maslov coorientation of the sphere of folds
$S$.
Suppose there is another function
$g_{{\alpha}^\prime , {\beta}^\prime , Q^\prime}$
of this type  such that the sign of
$\alpha^\prime ({\bf r})$
along
$S$
is the same as for
$\alpha ({\bf r})$.
Then it is not hard to show that there exists a symplectomorphism of
$\widetilde U$
preserving the fibers of the Lagrangian projection and mapping
the Lagrangian submanifold generated by
$g_{{\alpha}, {\beta}, Q}$
into the Lagrangian submanifold generated by
$g_{{\alpha}^\prime , {\beta}^\prime , Q^\prime}$.
Thus we have reduced 
$L$
and
$\cal F$
in a neighborhood of
$S$
to the standard form.
Therefore we can construct the necessary map
$F$
near each of the spheres
$S_1$
and
$S_2$.

Now we need to extend 
$F$
to a neighborhood of the rest of
$B$.
We can easily extend 
$F$
to a neighborhood of the whole 
$B$
as a
{\it
diffeomorphism
}
$\tilde F$
which maps 
${\cal F}$
into
${\cal F}^\prime$,
$L$
into
$L^\prime$
and which identifies the symplectic form 
$\omega$
on
$T^\ast L$
into the symplectic form 
$\omega^\prime$
on
$T^\ast L^\prime$
{\it
(only) along the zero-section
$L\subset T^\ast L$.
}

The form
$\omega_1 = {\tilde F}^\ast \omega^\prime$
is a symplectic form on a neighborhood 
$U$
of
$B$
in
$T^\ast L$
which already coincides with
$\omega$
along a subset
$W\subset U$
which is a union
of the zero-section
$L$
and of some neighborhoods of
$S_1$
and
$S_2$
in
$U$.
It suffices to find a symplectomorphism of a possibly smaller neighborhood 
$U_1\subset U$
of
$B$
which is identical on
$W$,
which preserves the leaves of the foliation
$\cal F$
and the zero-section
$L$
and which maps the form
$\omega$
on
${\widetilde U}$
into the form
$\omega_1$.
We find such a symplectomorphism using  Moser's argument.
Namely, let us consider a linear homotopy
$\omega_\tau = \tau\omega + (1-\tau) \omega_1$,
$0\leq\tau\leq 1$.
In a neighborhood of
$L$
this is a family of symplectic forms.
Observe that the leaves of  
${\cal F}$
are Lagrangian submanifolds with respect to
${\tilde F}^\ast \omega^\prime$
and hence with respect to any
$\omega_\tau$.
Observe also that we only need to construct the symplectomorphism near those points of
$L$
where
$\cal F$
is transversal to
$L$.
All this implies that we can solve the corresponding homological equation and find a (time-dependent) 
Hamiltonian vector field 
$X_t$
on a neighborhood
$U_1\subset U$
of
$B$
so that 

\smallskip
\noindent
(i) 
$X_t$
is zero at the points of
$W$;

\noindent
(ii) 
$X_t$
is tangent to
${\cal F}$;

\noindent
(iii)
the Hamiltonian flow 
$h_\tau$
of
$X$
is defined in a neighborhood of
$L\subset U_1$ 
and satisfies the property
$(h_\tau)^\ast \omega = \omega_\tau$.
\smallskip

The Hamiltonian symplectomorphism
$h_1$
combined with the diffeomorphism
$\tilde F$
gives the needed symplectomorphism between
a neighborhood of 
$B$
in
$T^\ast L$
and a neighborhood of
$B^\prime$
in
$T^\ast L^\prime$. 
This finishes the proof of the proposition.
\b

\end{document}